\documentclass{gtart_h}  

\def\ifplaintex{\expandafter\ifx\csname documentclass\endcsname\relax}

\def\ifplaintex{\expandafter\ifx\csname documentclass\endcsname\relax}

\ifplaintex 
\else
\fi

\expandafter\ifx\csname epsfbox\endcsname\relax\input epsf\fi

\def\gt{{\mathsurround=0pt\it $\cal G\mskip-2mu$eometry \&\ 
$\cal T\!\!$opology}}        

\def\gtp{{\mathsurround=0pt\it $\cal G\mskip-2mu$eometry \&\ 
$\cal T\!\!$opology $\cal P\!$ublications}}  

\def\lognumber#1{\def\thelognumber{#1}}
\def\volumenumber#1{\def\thevolumenumber{#1}}
\def\papernumber#1{\def\thepapernumber{#1}}
\def\volumeyear#1{\def\thevolumeyear{#1}}

\def\pagenumbers#1#2{\def\startpage{#1}\def\finishpage{#2}}
\def\published#1{\def\publishdate{#1}}
\def\proposed#1{\def\theproposer{#1}}
\def\seconded#1{\def\theseconders{#1}}
\def\received#1{\def\receiveddate{#1}}
\def\revised#1{\def\reviseddate{#1}}
\def\accepted#1{\def\accepteddate{#1}}

\def\coverauthors#1{\def\thecoverauthors{#1}}
\def\asciiauthors#1{\def\theasciiauthors{#1}}
\def\asciiaddress#1{\def\theasciiaddress{#1}}
\def\asciiemail#1{\def\theasciiemail{#1}}
\def\asciiurl#1{\def\theasciiurl{#1}}
\long\def\asciiabstract#1{\long\def\theasciiabstract{#1}}

\let\\\par\let\thelognumber\relax
\let\thevolumenumber\relax\let\thepapernumber\relax
\let\thevolumeyear\relax\let\thesamplenumber\relax\let\startpage\relax
\let\finishpage\relax\let\publishdate\relax\let\receiveddate\relax
\let\reviseddate\relax\let\accepteddate\relax\let\theasciititle\relax
\let\theasciiauthors\relax\let\theasciiaddress\relax
\let\theasciiabstract\relax
\let\theasciiemail\relax\let\theshortauthors\relax\let\theshorttitle\relax
\let\thecoverauthors\relax\let\theasciiurl\relax

\long\def\maketitlep{   

\count0=\startpage

\gt\hfill      
\hbox to 77pt{\vbox to 0pt{\vglue -15pt\epsfbox{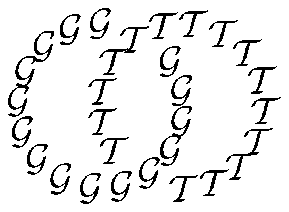}\vss}\hss}
\break
{\small\ifx\thesamplenumber\relax 
Volume \else Sample
\fi\thevolumenumber\ (\thevolumeyear)
\startpage--\finishpage\nl
Published: \publishdate}
\vglue 0.5truein plus 0.4fil minus 0.1truein

{\parskip=0pt\leftskip 0pt plus 1fil\def\\{\par\smallskip}{\ifplaintex\large
\else\Large\fi\bf\thetitle}\par\medskip}   

\vglue 0pt plus 0.1fil 

{\parskip=0pt\leftskip 0pt plus 1fil\def\\{\par}{\sc\theauthors}
\par\medskip}

\vglue 0pt plus 0.1fil 

{\small\parskip=0pt\let\newline\\
{\leftskip 0pt plus 1fil\def\\{\par}{\sl\theaddress}\par}
\expandafter\ifx\theemail\relax    
\relax\else\vglue 5pt plus 0.02fil minus 2pt\def\\{\stdspace{\rm 
and}\stdspace} 
\cl{Email:\stdspace\tt\theemail}\fi
\ifx\theurl\relax                  
\relax\else\vglue 5pt plus 0.02fil minus 2pt\def\\{\stdspace{\rm 
and}\stdspace}
\cl{URL:\stdspace\tt\theurl}\fi\par}

\vglue 7pt plus 0.3fil minus 3pt

{\bf Abstract}
\vglue 5pt plus 0.1fil minus 2pt

\theabstract

\vglue 7pt plus 0.3fil minus 3pt

{\bf AMS Classification numbers}\quad Primary:\quad \theprimaryclass

Secondary:\quad \thesecondaryclass

\vglue 5pt plus 0.3fil minus 2pt

{\bf Keywords:}\quad \thekeywords

\vglue 10pt plus 0.5fil minus 5pt

{\small  Proposed: \theproposer\hfill Received: \receiveddate\nl
Seconded: \theseconders\hfill 
\ifx\reviseddate\relax                         
Accepted: \accepteddate                        
\else
Revised: \reviseddate                          
\fi}
\eject
}       

\font\phead=cmsl9 scaled 950
\font=cmsl9 scaled 1050
\font\pnum=cmbx10 scaled 913
\font\lnum=cmbx10 
\font\pfoot=cmsl9 scaled 950
\font=cmsl9 scaled 1050
\ifplaintex
\headline{\vbox to 0pt{\vskip -4.5mm\line{\small\phead\ifnum
\count0=\startpage ISSN 1364-0380 (on line)
1465-3060 (printed) \hfill {\pnum\folio}\else\ifodd\count0\def\\{ }%
\ifx\theshorttitle\relax\thetitle\else\theshorttitle\fi\hfill{\pnum\folio}
\else\def\\{ and }{\pnum\folio}\hfill\ifx\theshortauthors\relax\theauthors
\else\theshortauthors\fi\fi\fi}\vss}}
\footline{\vbox to 0pt{\vglue 0mm\line{\small\pfoot\ifnum\count0=\startpage
\copyright\ \gtp\hfill\else
\gt, Volume \thevolumenumber\ (\thevolumeyear)\hfill\fi}\vss
}}
\else
\makeatletter
\def\@oddhead{{\small\lhead\ifnum\count0=\startpage ISSN 1364-0380 (on line)
1465-3060 (printed) \hfill {\lnum\number\count0}\else\ifodd\count0
\def\\{ }\ifx\theshorttitle\relax \thetitle \else\theshorttitle\fi\hfill
{\lnum\number\count0}\else\def\\{ and }{\lnum\number\count0}
\hfill\ifx\theshortauthors\relax 
\theauthors\else\theshortauthors\fi\fi\fi}}\def\@evenhead{\@oddhead}
\def\@oddfoot{\small\lfoot\ifnum\count0=\startpage\copyright\ \gtp\hfill\else
\gt, Volume \thevolumenumber\ (\thevolumeyear)\hfill\fi}
\def\@evenfoot{\@oddfoot}
\makeatother
\fi

\newwrite\gtoutfile
\long\gdef\makeheadfile{  
{\def\\{, }\def\s{ }
\immediate\openout\gtoutfile head.xxx
\immediate\write\gtoutfile{Proxy-for: \ifx\theasciiauthors\relax
\theauthors\else\theasciiauthors\fi\s<\ifx\theasciiemail\relax\theemail\else\theasciiemail\fi>}
\immediate\write\gtoutfile{\noexpand\\}
\immediate\write\gtoutfile{Authors: \ifx\theasciiauthors\relax
\theauthors\else\theasciiauthors\fi}
{\def\\{ }\immediate\write\gtoutfile{Title: \ifx\theasciititle\relax
\thetitle\else\theasciititle\fi}}
\immediate\write\gtoutfile{Subj-class: GT or SG or MG etc}
\immediate\write\gtoutfile{MSC-class: \theprimaryclass\ifx\thesecondaryclass\relax\else, \thesecondaryclass\fi}
\immediate\write\gtoutfile{Journal-ref: Geom. Topol. \thevolumenumber
(\thevolumeyear) \startpage-\finishpage}
\immediate\write\gtoutfile{Comments: Published by Geometry and Topology at}
\immediate\write\gtoutfile{\s\s http://www.maths.warwick.ac.uk/gt/GTVol\thevolumenumber/paper\thepapernumber.abs.html}
\immediate\write\gtoutfile{\noexpand\\}
\immediate\write\gtoutfile{}
\ifx\theasciiabstract\relax
\immediate\write\gtoutfile{\theabstract}\else
\immediate\write\gtoutfile{\theasciiabstract}\fi
\immediate\write\gtoutfile{}
\immediate\write\gtoutfile{\noexpand\\}
\immediate\write\gtoutfile{}
\immediate\closeout\gtoutfile}}  

\def\maketitlepage{\maketitlep\makeheadfile}
\let\maketitle\maketitlepage

\lognumber{450}
\received{30 April 2004}
\volumenumber{9}\papernumber{3}\volumeyear{2005}
\pagenumbers{121}{162}   
\revised{22 December 2004}
\published{28 December 2004}
\accepted{27 December 2004}
\proposed{Ralph Cohen}
\seconded{Leonid Polterovich, Frances Kirwan}

\usepackage{amsmath,amssymb,amscd,euscript}
\usepackage[all]{xy}

\def\S{Section }

\newcommand{\T}{{\mathbb T}}

\newcommand{\Ta}{{\widetilde a}}

\newcommand{\B}[1]{{\mathbb #1}}

\newtheorem{theorem}{Theorem}[section]

\newtheorem{cor}[theorem]{Corollary}
\newtheorem{lemma}[theorem]{Lemma}
\newtheorem{proposition}[theorem]{Proposition}
\newtheorem{prop}[theorem]{Proposition}

\newtheorem{defn}[theorem]{Definition}
\newtheorem{example}[theorem]{Example}
\newtheorem{remark}[theorem]{Remark}
\newtheorem{rmk}[theorem]{Remark}

\numberwithin{figure}{section}
\numberwithin{equation}{section}
\numberwithin{table}{section}

\newcommand{\p}{{\partial}}
\newcommand{\al}{{\alpha}}

\newcommand{\be}{{\beta}}

\newcommand{\Om}{{\Omega}}
\newcommand{\om}{{\omega}}
\newcommand{\eps}{{\varepsilon}}

\newcommand{\ga}{{\gamma}}
\newcommand{\Ga}{{\Gamma}}
\newcommand{\io}{{\iota}}
\newcommand{\ka}{{\kappa}}
\newcommand{\la}{{\lambda}}
\newcommand{\La}{{\Lambda}}
\newcommand{\si}{{\sigma}}
\newcommand{\Si}{{\Sigma}}
\newcommand{\cp}{{\B C\B P}}
\newcommand{\Kk}{{\mathcal K}}
\newcommand{\Ll}{{\mathcal L}}
\newcommand{\Aa}{{\mathcal A}}
\newcommand{\Ii}{{\mathcal I}}
\newcommand{\tr}{{\rm tr\,}}
\newcommand{\ov}{\overline}
\newcommand{\THha}{{\widetilde{\Hh}\,\!^a_0}}

\newcommand{\del}{\partial}
\newcommand{\Rar}{\longrightarrow}

\newcommand{\map}[1]{\stackrel {#1}\longrightarrow}

\renewcommand{\Tilde}{\widetilde}

\newcommand{\Cc}{{\mathcal C}}
\newcommand{\Gg}{{\mathcal G}}
\newcommand{\Dd}{{\mathcal D}}
\newcommand{\Mo}{(M,\omega)}
\newcommand{\Ma}{(M,a)}
\newcommand\Hom{\operatorname{Hom}}

\newcommand\Flux{\operatorname{Flux}}
\newcommand\Symp{\operatorname{Symp}}
\newcommand\Ham{\operatorname{Ham}}

\newcommand\im{\operatorname{Im}}

\newcommand\Diff{\operatorname{Diff}}
\newcommand\Lie{\operatorname{Lie}}
\newcommand{\Hh}{{\mathcal H}}
\newcommand{\Q}{{\mathbb Q}}
\newcommand{\ev}{{\rm ev}}

\newcommand{\R}{{\mathbb R}}
\newcommand{\Z}{{\mathbb Z}}
\newcommand{\oq}{{\otimes \Q}}

\newcommand{\Tf}{{\widetilde f}}

\newcommand{\ms}{{\medskip}}

\begin{document}
    
\title{Homotopy properties of Hamiltonian group actions}
\author{Jarek K\c edra\\Dusa McDuff}
\coverauthors{Jarek K\noexpand\c edra\\Dusa McDuff}
\asciiauthors{Jarek Kedra\\Dusa McDuff}

\address{Institute of Mathematics US, Wielkopolska 15, 70-451 
Szczecin, Poland\\{\rm and}\qua  Mathematical Institute Polish 
Academy of Sciences\\\'Sniadeckich 8, 00-950 Warszawa, Poland
\\{\rm and}\\Department of Mathematics,
 Stony Brook University\\Stony Brook, 
NY 11794-3651, USA\\\smallskip\\{\rm Email:\qua}{\tt\mailto{kedra@univ.szczecin.pl}{\rm\qua
and\qua}\mailto{dusa@math.sunysb.edu}}\\\smallskip\\\tt\url{http://www.univ.szczecin.pl/~kedra}, \url{http://www.math.sunysb.edu/~dusa}}

\asciiaddress{Institute of Mathematics US, Wielkopolska 15, 70-451 
Szczecin, Poland\\and  Mathematical Institute Polish 
Academy of Sciences\\Sniadeckich 8, 00-950 Warszawa, Poland
\\and\\Department of Mathematics,
 Stony Brook University\\Stony Brook, NY 11794-3651, USA}

\asciiemail{kedra@univ.szczecin.pl, dusa@math.sunysb.edu}

\asciiurl{http://www.univ.szczecin.pl/ kedra,
http://www.math.sunysb.edu/ dusa}


\begin{abstract}  Consider a Hamiltonian action of
a compact Lie group $G$ on a compact symplectic manifold $\Mo$ and 
let $\Gg$ be a subgroup of the diffeomorphism group $\Diff M$.
We develop techniques to decide when
the maps on rational homotopy and rational homology 
induced by the classifying map $BG\to B\Gg$ are injective.
For example, we extend Reznikov's result for
complex projective space $\cp^n$ 
to show that both  in this case and the case of
generalized flag manifolds the 
natural map $H_*(BSU(n+1))\to H_*(B\Gg)$ is injective,
where $\Gg$ denotes the group of all diffeomorphisms 
that act trivially on cohomology.
We also show that if $\la$ is a Hamiltonian circle action 
that contracts in $\Gg: = \Ham\Mo$ then there is an associated nonzero 
element in $\pi_3(\Gg)$ that deloops to a nonzero element of 
$H_4(B\Gg)$.  This result 
(as well as many others)
extends to c-symplectic 
manifolds $\Ma$, ie, $2n$--manifolds with a class $a\in 
H^2(M)$ such that $a^n\ne0$. 
The proofs are based on 
calculations of certain characteristic classes
and elementary homotopy theory. 
\end{abstract}

\asciiabstract{%
Consider a Hamiltonian action of a compact Lie group H on a compact
symplectic manifold (M,w) and let G be a subgroup of the
diffeomorphism group Diff(M).  We develop techniques to decide when
the maps on rational homotopy and rational homology induced by the
classifying map BH --> BG are injective.  For example, we extend
Reznikov's result for complex projective space CP^n to show that both
in this case and the case of generalized flag manifolds the natural
map H_*(BSU(n+1)) --> H_*(BG) is injective, where G denotes the group
of all diffeomorphisms that act trivially on cohomology.  We also show
that if lambda is a Hamiltonian circle action that contracts in G =
Ham(M,w) then there is an associated nonzero element in pi_3(G) that
deloops to a nonzero element of H_4(BG).  This result (as well as many
others) extends to c-symplectic manifolds (M,a), ie, 2n-manifolds with
a class a in H^2(M) such that a^n is nonzero.  The proofs are based on
calculations of certain characteristic classes and elementary homotopy
theory.}

\maketitle

\section{Introduction}\label{S:intro}

\subsection{Overview of results}

This paper studies  the homotopy type of the group
$\Symp\Mo$   of symplectomorphisms of
a closed symplectic manifold $\Mo$ onto itself.  It was noticed from 
the beginning of  the modern development of symplectic topology 
that many basic results  can be expressed in terms of 
properties of this group.  For instance, one way to express 
symplectic rigidity is to observe that 
the group $\Symp\Mo$ is closed in the $C^0$ (or uniform) topology 
on the diffeomorphism group $\Diff (M)$, rather than just 
in the $C^1$--topology  as one would expect.

In dimension two, it follows from Moser's theorem that
the group of 
oriented area preserving
(or symplectic) diffeomorphisms of any surface $(\Si, \om)$ 
is homotopy equivalent to the full group of orientation 
preserving diffeomorphisms $\Diff^+(\Si)$.
Moreover, its homotopy type is well known.  For example, in 
the case of the $2$--sphere, Smale showed that 
$\Symp(S^2, \om)\simeq \Diff^+(S^2)$ is homotopy equivalent to
the group  $SO(3)$ of rotations.

In dimension four, the two groups $\Symp\Mo$ and $\Diff^+(M)$
are very different. Very little is known about the homotopy groups of 
$\Diff^+(M)$ (even for the $4$--sphere) while in certain cases it
has  turned out that $\Symp\Mo$ is accessible.  Using  the method of
pseudoholomorphic curves, Gromov~\cite{G} proved that the 
identity component of the group 
$\Symp(S^2\times S^2,\om_1)$ is homotopy equivalent to 
to $SO(3)\times SO(3)$.
Here $\om_{\la} $ denotes the product symplectic form in which the first sphere
$S^1\times \{p\}$ has area equal to $\la$ and the second $\{p\}\times S^2$ 
has area equal to $1$. Hence in the above situation both spheres have equal size.
Similarly, the symplectomorphism group of $\cp^2$ with its 
standard symplectic form deformation retracts to the group of linear isometries $PSU(3)$.  

Abreu, Anjos and McDuff~\cite{AM,An} extended Gromov's work, studying the family of groups 
$\Symp(S^2\times S^2,\om_\la)$ for $\la>1$.  They found that 
although these groups are not homotopy equivalent to any
compact Lie group, their rational homotopy is
detected by actions of compact Lie groups. More precisely,
the rational homotopy of $\Symp(S^2\times S^2,\om_{\la})$
is nontrivial in degrees $1,3$ and $4k$ where $k$
is the largest integer $<\la$.
The generator in degree one is represented by a circle
action, the two generators in degree three are represented
by actions of $SO(3)$ and the last generator is a certain higher order
Samelson product of the previous ones.\footnote
{Recall that if $G$ is a topological group then 
the Samelson product
$\bigl<,\bigr >\co \pi_k(G)\times \pi_m(G)\to \pi_{k+m}(G)$
is defined by\nl
\cl{$\bigl<\al,\be\bigr> \co S^{k+m}=S^k \times S^m /S^k\vee S^m
\to G.\quad [s,t]\mapsto\al(s)\be(t)\al(s)^{-1}\be(t)^{-1}.$}}
This set of results
culminated in a theorem of Anjos--Granja \cite{AG} stating
that when  $1<\la\leq 2$ the group $\Symp(S^2\times S^2,\om_{\la})$ is homotopy equivalent
to the topological amalgamated product 
$$
(SO(3)\times SO(3)) *_{SO(3)} (SO(3)\times S^1),
$$
where $SO(3)$ maps to the diagonal subgroup of $SO(3)\times SO(3)$ and 
is included as the first factor in $SO(3)\times S^1$.

In all these cases the rational homotopy of the symplectomorphism 
group is determined by some compact subgroups arising from
Lie group actions that preserve 
a  K\"ahler metric on the underlying manifold.
As a first step towards understanding what happens in
higher dimensions, 
 we look at manifolds that admit an effective  
symplectic action by a compact
connected
 Lie group $G$ and see how much 
of the (rational) homotopy  of $G$
remains visible in $\Symp\Mo$, either directly in the sense that 
the induced map on $\pi_*\oq$ is injective, or indirectly 
in the sense that there is some associated nontrivial 
element in $\pi_*(\Symp\Mo)$ (such as a secondary Samelson product).  
Reznikov~\cite{Rez} proved some 
initial results in this direction.
By defining and calculating some new characteristic classes, 
he showed that any 
nontrivial homomorphism  $SU(2)\to \Symp\Mo$ induces an 
injection on $\pi_3$. He also showed that
the canonical map of $SU(n+1)$
into $\Symp(\cp^n)$ induces an injection on rational homotopy.
In fact, he proved the sharper result that for $i>1$ the Chern 
classes $c_i\in H^{2i}(BSU(n+1);\R)$ have natural extensions to classes in 
$H^{2i}(B\Symp(\cp^n);\R)$.  Thus the induced map
$$
H_*(BSU(n+1);\Q)\to H_*(B\Symp(\cp^n);\Q) 
$$
is injective.\footnote
{
The injectivity of $r_*\co \pi_*(G)\oq\to \pi_*(\Gg)\oq$ is equivalent  to  that of the induced map on classifying spaces
$R_*\co \pi_*(BG)\oq\to \pi_*(B\Gg)\oq$. However this does not mean that
the induced map on rational homology
$H_*(BG;\Q)\to H_*(B\Gg;\Q)$ must be injective, because  the rational Hurewicz map
$h\co \pi_*(B\Gg)\oq\to H_*(B\Gg)\oq$ is not always injective when 
$\Gg=\Symp\Mo$, though it is for connected Lie groups. Equivalently,
 the rational Whitehead product need not vanish on $B\Gg$:  see \S\ref{S:product}.}
Similar results were obtained for certain toric manifolds by 
Janusz\-kiewicz--K\c edra in~\cite{JK}.

In \S\ref{S:ccl} of this paper we show that Reznikov's characteristic 
classes are closely related to those used by Januszkiewicz--K\c edra
and to the $\ka$--classes of Miller--Morita--Mumford.   
In particular they may be defined without using any  
geometry, and so they extend to {\em cohomologically symplectic (or c-symplectic)}
manifolds $(M,a)$, ie, to pairs consisting of
a closed oriented $2n$--manifold $M$ together with a
cohomology class $a\in H^2(M;\R)$ such that $a^n>0$.
However, if $H^1(M;\R)\ne 0$, these classes do not live on 
the analog of the  full symplectomorphism group but rather 
on the analog of its Hamiltonian\footnote
{
This consists of the time-$1$ maps $\phi_1^H$ of Hamiltonian flows 
$\phi_t^H, t\in [0,1]$. These paths are generated by time 
dependent Hamiltonian functions $H_t\co M\to \R$ via the recipe $\om(\dot\phi_t^H, \cdot)= 
dH_t(\cdot)$. If $H^1(M;\R) = 0$, then $\Ham\Mo$ is just the 
identity component $\Symp_0\Mo$ of the  symplectomorphism group.  
In general, it is the kernel of the flux homomorphism~(\ref{eq:flux}): 
for further background information and references, 
see for example McDuff--Salamon~\cite{MS}.} 
subgroup $\Ham\Mo$. If $H^1(M;\R)= 0$ they can be extended significantly further, 
see Remark~\ref{rmk:s0} and \S\ref{S:ccl}.

We extend Reznikov's work in two directions, first by looking at 
homogeneous spaces other than $\cp^n$, and second by considering 
homotopy groups in dimensions other than three.  The first is fairly
straightforward; for example we prove in
Proposition~\ref{prop:flags2} that the
action of $SU(n)$ on generalized flag manifolds induces an
injection of rational homology at the classifying space level.
However, the second  is more delicate because it
is not true that every inclusion $U(k)\to \Symp\Mo$
induces an injection on $\pi_i\oq$ for $i\ne 3$.
As shown in more detail in \S\ref{ss:ex}, our counterexamples are based on the existence of homomorphisms
$U(m)\to U(n)$ that kill homotopy.
For instance the homomorphism 
$$
U(m)\to U(m)\times U(m)\subset U(2m),\quad
A\mapsto (A,\ov A)
$$
(where $\ov A$ is
complex conjugation) kills homotopy in dimension $4k+1$.
One cannot remedy this by considering maximal compact subgroups
of $\Symp\Mo$:
McDuff--Tolman~\cite{MT2} construct a toric
K\"ahler structure on the product $\cp^1\times \cp^2$ whose isometry group $G$ is maximal in $\Symp\Mo$ but is such that $\pi_1(G)\oq$
does not inject.
Therefore the homotopy of $G$ need not be directly visible in  $\Symp$.
However, we shall see in Theorem~\ref{thm:s0} that in the case of a null homotopic circle action one 
can use the null homotopy itself to define 
an element $\rho\in\pi_3(\Symp\Mo)$ that never vanishes.
This is 
the simplest case of a general construction that is developed in \S\ref{S:product}.
In \S\ref{sec:eval} we discuss further properties of this element $\rho$, 
looking in particular at its image in $\pi_3(M)$ under the pointwise evaluation 
map $\Symp\Mo\to M$.

This paper uses elementary methods from algebraic topology.  
For recent results on the homotopy of $\Symp\Mo$ that use deeper, more analytic 
techniques, readers might consult McDuff~\cite{Mcox} or Seidel~\cite{SeiD}. 
In the rest of this introduction we state our main results in more detail. 
We denote by $H_*(M), H^*(M)$ (co)homology with real or rational 
(rather than integral) coefficients.

\subsection{Circle actions}

Let $\Gg$ be a topological group and $\la\co S^1\to \Gg$
a nonconstant homomorphism that represents the zero element in 
$\pi_1(\Gg)$ and so extends to a map $\Tilde\la\co D^2\to \Gg$.  
(For short, we say that $S^1$ is {\em inessential} in $\Gg$.)
Define $\rho \in \pi_3(\Gg)$  
by setting
\begin{eqnarray}\label{eq:rho3}
&&\rho\;\co \;S^3 = 
\bigl(D^2\times S^1\bigr)/\bigl((D^2\times \{1\})\vee 
(\p D^2\times S^1)\bigr)\to \Gg,  
\\\notag
&&\qquad\qquad\qquad\quad
(z,t) \mapsto \left<\Tilde\la(z),\la(t)\right>,
\end{eqnarray}
where the bracket $\left<\phi,\psi\right>$ represents the commutator
$\phi\psi\phi^{-1}\psi^{-1}$.  
Observe that the map $D^2\times S^1\to \Gg$ descends to $S^3$ precisely because 
it contracts the boundary $\p D^2\times S^1$ to a point.  
Thus it is crucial here that $G=S^1$ is abelian
and that $\la$ is a homomorphism.

\begin{theorem}\label{thm:s0}  Let $(M, \om)$ be a closed symplectic 
manifold of dimension $2n$ and set $\Gg: = \Ham(M,\om)$.
Let $\la \co S^1 \to \Gg$ be a nontrivial homomorphism that  is 
inessential in $\Gg$.
Then the element $\rho\in \pi_3(\Gg)$ defined above is independent of the choice of extension $\Tilde\la$ and has infinite order. 
Moreover, the corresponding
element
$\overline{\rho}\in \pi_4(B\Gg)$  has
nonzero image in $H_4(B\Gg;\B Q)$.
\end{theorem}

\begin{rmk}\label{rmk:s0}\rm 
(i)\qua To see that $\rho$ is independent of the choice of 
$\Tilde\la$, observe that any two extensions differ by an element $\be\in
\pi_2(\Gg)$.  This changes  $\rho$ 
by the Samelson product $\langle\be,\la\rangle$, which vanishes in $\Gg$ since $\la = 0$ in $\pi_1(\Gg)$.\smallskip

(ii)\qua
If $\la$ is a smooth inessential circle action on $M$ then one can
use formula~(\ref{eq:rho3}) to define an element
$\rho\in \pi_3(\Diff(M))$.  However, in this generality
we have no way of proving that $\rho\ne 0$.
\smallskip

(iii)\qua 
We detect the nontriviality of $\ov{\rho}$ by using 
the characteristic classes of Reznikov~\cite{Rez}
and Januszkiewicz--K\c edra~\cite{JK}.  These classes extend beyond 
the Hamiltonian group $\Ham\Mo$ to appropriate
topological monoids $\Hh$ that act on $M$.
For example,  we show in Corollary~\ref{cor:nonzero0}  that if 
$\la$ is a smooth action on a simply connected
c-symplectic manifold $\Ma$ 
that is inessential in the topological monoid $\Hh_a$ formed by all smooth homotopy 
equivalences $M\to M$ that 
fix the class $a$ then $\ov{\rho}$ has 
nonzero image in $H_4(B\Hh_a)$.
\smallskip

(iv)\qua 
Our result extends work by Reznikov in the following way.
Reznikov proved in~\cite{Rez}
that  any Hamiltonian action of $SU(2)$ induces a nonzero map
on $\pi_3\oq$.   Moreover,  
Lemma~\ref{le:un} below implies that the element $\rho$ created from 
any circle subgroup of $SU(2)$ lies in the image of $\pi_3(SU(2))$.
Thus the nontriviality statement in Theorem~\ref{thm:s0}
follows from Reznikov's work in the case when the 
circle $\la$ contracts in $\Ham\Mo$ by virtue of the fact that it is 
contained in a simply connected Lie subgroup $G$ of $\Ham\Mo$.  
\smallskip

(v)\qua 
Let $\ev_*\co \pi_*(\Gg)\to \pi_*(M)$ 
denote the map obtained by 
evaluating at the base point $p$.  By looking at the $SU(2)$--action on 
$M:=\cp^n$ for $n=1,2$ 
and using Lemma~\ref{le:un} as in (iv) above, one sees  
that the element $\ev_*(\rho)\in 
\pi_3(M)$ is sometimes zero and sometimes nonzero.  We show in 
Proposition~\ref{prop:al2} that if  $H^1(M) = 0$ then
$\ev_*(\rho) \ne 0$  only if there is a nonzero
quadratic relation $\sum_{ij} c_ic_j=0$ among the classes $c_i\in
H^2(M)$. In the symplectic case this is no surprise since it follows 
from the work of Lalonde-McDuff~\cite{LM} that the map 
$$
h\circ\ev_*\co  \pi_*(\Ham\Mo)\to H_*(M;\Q)
$$
is zero, where $h$ denotes the Hurewicz homomorphism.  Hence, by minimal model theory, any element in the image 
of $\ev_*$ in this low degree 
must give rise to a relation in $H^*(M)$.
However, the arguments in~\cite{LM} do not apply in the 
c-symplectic case.
\qed
\end{rmk}

A Hamiltonian $S^1$--action always has fixed points $p$
and one can consider its image in the subgroup 
$\Diff_p$ of  diffeomorphisms that fix $p$.  When $H^1(M)\ne 0$ 
we shall need to consider the corresponding subgroup $\Gg_p$ of 
$\Gg: = \Ham\Mo$
which is the fiber of the evaluation map $\Gg\to 
M$ at $p$. The following result is proved in \S\ref{ss:circ}.

\begin{lemma}\label{le:p}  Suppose that $\la$ is a  
Hamiltonian circle action with 
moment map $H\co M\to \R$, and set 
$\Gg: = \Ham\Mo$.
 Then $\la$ is essential in 
    $\Gg_p$ for every fixed point $p$ such that
    $\int_M (H-H(p))\om^n \ne 0$.  
    In particular, $\la$ is essential in $\Gg_p$ if $H$ assumes its 
    maximum or minimum at $p$.
    \end{lemma}

A similar statement holds in the c-Hamiltonian case.

\subsection{Higher homotopy groups}

The proof of Theorem~\ref{thm:s0} is based on
a general construction of a secondary product $\{\cdot,\cdot\}$ on the
homotopy of topological groups. Given a point $p\in M$ we shall write\footnote
{Our conventions are that $G$ denotes a compact Lie group while $\Gg$ (usually) denotes an infinite dimensional group such as $\Symp$ or $\Diff$.}
$G_p, \Gg_p$ for the subgroups of $G,\Gg$ that fix $p$.  

Let $r\co G\to \Gg$ be a continuous homomorphism.
Consider the induced map
$$
r_{k-1}\co \pi_{k-1}(G_p)\oq\to \pi_{k-1}(\Gg)\oq.
$$
In \S\ref{S:product} we construct an element 
$$
\{f,f'\}\in \pi_{k+m-1}(\Gg)/({\rm im\,}r_{k+m-1})
$$
for each pair $f\in \ker r_{k-1}, f'\in \pi_{m-1}(G_p)$ 
that is well defined modulo Samelson products of the form 
$\langle\be,r\circ f'\rangle$, 
$\be\in \pi_k(\Gg)$.  The most important case is when 
$f$ is given by an inessential circle action.  Lemma~\ref{le:un} 
states that rational homotopy of the group $SU(n)$ is generated by products of this form.
In \S3 we use characteristic classes
to calculate $\{f,f'\}$ 
in various other cases, for example when $f=f'$ is given by 
an inessential  Hamiltonian circle action.

The next proposition follows by combining Lemma~\ref{le:un} with the proof of Theorem~\ref{thm:s0}.
We shall state it in its most general form, ie, for a   
c-symplectic manifold $\Ma$ and for the largest possible topological monoid $\Hh$.  If $H^1(M;\R)= 0$, we may take $\Hh$ to be the monoid
$\Hh_a$ mentioned in Remark~\ref{rmk:s0}(iii).  However, 
as we explain in more detail in \S\ref{S:ccl}, if  
$H^1(M;\R)\ne 0$ then we must work in a context in which the $a$--Flux homomorphism $\Flux^a$ vanishes. Let $\Hh_0$ denote the identity component of $\Hh_a$, ie,
 the monoid formed by all smooth maps $M\to M$ that are homotopic to the identity.  Then $\Flux^a\co  \pi_1(\Hh_0)\to H^1(M;\R)$
is given by
\begin{equation}\label{eq:aflux}
  \Flux^a(\la)(\ga) = \langle a, \tr_\la(\ga)\rangle\qquad \mbox{for }
 \ga\in H_1(M),
\end{equation}
 where $\tr_\la(\ga)\in H_2(M)$ is represented by 
 $T^2\to
 M,(s,t)\mapsto \la(s)\bigl(\ga(t)\bigr)$.  We shall denote by
 $\THha$ the cover of
$\Hh_0$ determined by the $a$--Flux homomorphism.  Since $\Flux^a$ is a homomorphism,  $\THha$ is a monoid.
Note also that in the symplectic case, $\Flux^a$ vanishes on all loops in $\Ham\Mo$.  Hence 
the map $g\mapsto (g,[g_t])$ (where $\{g_t\}_{t\in[0,1]}$ is any path in $\Ham$ from the identity to $g$) defines an inclusion $\Ham\Mo\to \THha$.

\begin{prop}\label{prop:nonzero}  Suppose that $G: = SU(2)$ acts 
smoothly and with finite kernel
on a c-symplectic manifold $(M,a)$. 
Then the induced map
$$
\pi_{4}(BG)\to H_{4}(B \THha)
$$
is injective. Further, if $H^1(M)=0$ 
then the image of $\pi_{4}(BG)$
in $H_{4}(B\Hh_a)$ is nonzero.\end{prop}

If $(M,a)$ supports an action of $G: = SU(\ell)$ with $\ell>2$
then one can also try 
to understand the maps $\pi_{2i}(BG)\oq\to H_{2i}(B\Gg)$ for $2< i 
\le \ell$ and appropriate $\Gg$.  Examples~\ref{ex:ex1} and~\ref{ex:ex2} show that these maps need not be injective.
It is still possible that the injectivity of 
$R_*\co  \pi_{2i}(BG)\oq\to \pi_{2i}(B\Gg)\oq$ implies that of 
$h\circ R_*\co   \pi_{2i}(BG)\oq\to H_{2i}(B\Gg)$, where
$h$ denotes the rational Hurewicz map $\pi_{2i}(BG)\oq\to H_{2i}(BG)$.\footnote{
  Note that if $\Gg$ were a compact Lie group this would be obvious, since in this case 
all rational Whitehead products in $B\Gg$ vanish.  
However, it was shown in Abreu--McDuff~\cite{AM} that Whitehead products do not always vanish when $\Gg: = \Symp\Mo$.  Hence in the generality considered here we cannot assume that $h$ is injective.}

Here is one result in this direction. 
Again we denote  by
$\Hh_0: = (M^M)_0$, the space of 
(smooth) self maps of $M$ that are homotopic to the identity.
Since the point evaluation
map $\ev_*\co \pi_*(G)\to \pi_*(M)$ factors through $\pi_*(\Hh_0)$ 
the kernel of $r_*\co  \pi_*(G)\to \pi_*(\Hh_0)$ is contained in the 
kernel of $\ev_*$.
We now describe conditions under which the kernel of the corresponding map 
$$
h\circ R_*\co  \pi_*(BG)\oq\to H_*(B\Hh_0)
$$
on the classifying space level is contained in ${\rm ker}\,\ev_*$.

We say that $M$ satisfies the
 {\em c-splitting condition}
  for a topological monoid $\Hh$ that acts on $M$ 
 if the associated 
 (Hurewicz)
 fibration $M\to 
   M_\Hh\to B\Hh$ is c-split, that is, if $\pi_1(B\Hh)$ acts trivially 
   on $H^*(M)$ and the Leray--Serre spectral sequence 
   for rational cohomology degenerates at the $E_2$--term.  
For example, Blanchard showed in~\cite{Bl} that
a simply connected K\"ahler manifold satisfies the 
c-splitting condition for $\Hh_0$: see Lemma~\ref{le:Bl}.
Further, it was conjectured in~\cite{LM} that any symplectic manifold 
satisfies the 
c-splitting condition for $\Ham\Mo.$
Recall that $M$ is called nilpotent if 
$\pi_1(M)$ is nilpotent and its action on $\pi_*(M)\oq$ is nilpotent
\cite{TO}.
The following result is proved in~\S\ref{ss:split}.
For $\al\in \pi_{2k-1}G$ we denote by $\ov{\al}$ the corresponding 
element in $\pi_{2k}(BG)$.

 \begin{prop}\label{prop:eval}  Let $G$ be a compact, connected, 
      and simply connected
    Lie group that acts (with finite kernel) 
 on a nilpotent manifold $M$.  Suppose further 
 that $M$ satisfies
 the c-splitting condition with respect to a 
 monoid $\Hh$  containing $G$. 
 Then for all $\al\in \pi_*(G)$
 $$
 \ev_*(\al) \ne 0\mbox{ in } \pi_{2k-1}(M)\oq 
 \;\;\Rightarrow\;\; h\circ R_*(\ov{\al}) \ne 0 \mbox{ in } H_{2k}(B\Hh).
 $$
\end{prop}

In this proposition $M$ can be an arbitrary smooth manifold.  
However, if it is not c-symplectic there is no obvious way 
to define an interesting monoid $\Hh$  that satisfies the c-splitting condition.
Proposition~\ref{prop:eval2} is an alternative version 
with a weaker splitting condition that holds when
$M$ is  symplectic.

In the next results $\Hh_H$ denotes the submonoid of $\Hh$
that acts trivially on cohomology.

\begin{cor}\label{cor:cfl}  Consider the action of $G: = SU(\ell)$ 
    on the manifold 
$M: = G/T$ of complete flags, where $T$ is the 
maximal torus in $G$. 
Then the map $H_*(BG)\to H_*(B\Hh_H)$ is injective.
\end{cor}
\begin{proof}\,  It follows from Lemma~\ref{le:Bl} that $M$ satisfies 
the c-splitting condition for $\Hh_H$.  Further  $\ev_*\co \pi_{*}(G)\oq\to 
\pi_{*}(M)\oq$ is injective. Hence
$\pi_*(BG)\oq\to H_*(B\Hh_H)$ is injective by Proposition~\ref{prop:eval}.
Since $H^*(BG)$ is freely generated by the duals of the 
spherical classes, the induced map $R^*\co  H^*(B\Hh)\to H^*(BG)$
is surjective, and so the corresponding map on rational homology is injective.
\end{proof}

A general flag manifold can be written as 
$$
M(m_1,\dots,m_k): = U(\ell)/U(m_1)\times\cdots\times U(m_k),\qquad
m_1\ge\dots\ge m_k,
$$
where $\ell = \sum m_i$. 
In Proposition~\ref{prop:flags2} below we extend the 
result of Corollary~\ref{cor:cfl} to arbitrary flag 
manifolds.  However the proof is considerably more complicated, and uses the characteristic classes defined in \S\ref{S:ccl}.
As a warmup we give in \S\ref{S:ccl}
a new proof of the following extension of
Reznikov's result about projective space.

\begin{prop}\label{prop:cpn}
In the case of the action of $G = SU(n+1)$ on $M: = \cp^n$,
the induced map
$
R_*\co H_*(BG)\oq\to H_*(B\Hh_H)$ is injective.
\end{prop}

The last section  \S\ref{sec:eval}
investigates the image $\ev_*(\rho)$ of the element $\rho$ 
of Theorem~\ref{thm:s0} under the evaluation map, and shows its 
relation to certain Whitehead products in $\pi_*(M)$.
For example, in Proposition~\ref{prop:al2} we give necessary and 
sufficient conditions for
 $\ev_*(\rho)$ to be nonzero, while 
 Proposition~\ref{prop:wh2} gives conditions under 
which the map $\ev_*\co  \pi_5(\Gg)\to \pi_5(M)$ is nonzero.

\medskip

{\bf Acknowledgements}\qua   We warmly thank the referees for making 
several important suggestions that have helped to simplify and clarify
various arguments, in particular the proofs of 
Propositions~\ref{prop:repff'}, \ref{prop:ka} and \ref{prop:flags2}. 
The first author thanks Thomas Vogel for discussions.
   
The first author is a member of EDGE, Research Training 
Network HPRN-CT-2000-00101, supported by the European Human Potential 
Programme. He is also partly supported by the KBN grant 1P03A 023 27.
The second author is partly supported by the NSF grant DMS 0305939.


\section{A secondary product on homotopy groups of classifying spaces}
\label{S:product}

If $\Gg$ is a topological group the
Samelson product  $\left <f,f' \right >\in \pi_{i+j}(\Gg)$ of the
elements $f\in \pi_i(\Gg), f'\in \pi_j(\Gg)$ is given by the map
\begin{equation}\label{eq:sam}
S^{i+j} = S^i\times S^j/S^i\vee S^j \to \Gg\co \quad
(x,y)\mapsto \left <f(x),f'(y)\right>.
\end{equation}
Here $\left<a,b\right>$ denotes the commutator $aba^{-1}b^{-1}$,
and we take the base
point of $\Gg$ to be the identity element so that the commutator
vanishes on the wedge $S^i\vee S^j$.
The main fact we shall use about
this product is that (up to sign)
it is the desuspension 
of the Whitehead product.
Thus the obvious isomorphism $\pi_*(\Gg)\cong \pi_{*+1}(B\Gg)$
takes $\pm \left <f,f' \right >$ to the Whitehead product $[F,F']$ of
the images $F, F'$ of $f, f'$.  Another important fact is that
the Whitehead product vanishes rationally on $BG$ when
$G$ is a compact connected Lie group.  This follows from minimal model
theory: nonzero Whitehead products lie in the kernel of the Hurewicz
homomorphism and give rise to relations in rational
cohomology, but $H^*(BG;\Q)$ is a free algebra.
Hence the rational Samelson product also vanishes on Lie groups,
but need not vanish when $\Gg$ is a
symplectomorphism group.

All results in this  section are concerned with rational homotopy and homology.
One could therefore work in the rational homotopy category, in which 
case
the notation $S^k$ does not denote a sphere but rather
its image $S^k_\Q$ in this category. Alternatively, one can
take sufficiently high
multiples of all maps so that they are null homotopic rather than zero
in $\pi_*\oq$, and then work in usual category.  We adopt the latter
approach.

\subsection{A general construction}

The first main result of  this section is Proposition~\ref{prop:uniq}.
In Corollary~\ref{cor:rho}
we establish the uniqueness 
(but not the nontriviality) of the element $\rho$ of 
Theorem~\ref{thm:s0}.

Consider the sequence of fibrations
associated with this action:
$$
\CD
\Gg_p @>i>> \Gg@>ev >> M @>j>> B\Gg_p = M_\Gg @>\pi>> B\Gg,
\endCD
$$
where $M_\Gg$ denotes the total space of the universal $M$--bundle over 
$B\Gg$.

Let $f\co S^{k-1}\to \Gg_p$ and $f'\co S^{m-1}\to \Gg_p$ be maps satisfying 
the
following assumptions:\ms

{\bf (A1)}\qua The Samelson product is trivial, $\left <f,f' \right > = 
0$;

{\bf (A2)}\qua The map $i\circ f$ is null homotopic in $\Gg$.\ms

Equivalently, for the corresponding maps to the classifying space,
$F\co S^k\to B\Gg_p$ and $F'\co S^m\to B\Gg_p$ we have\ms

{\bf (B1)}\qua The Whitehead product is trivial, $\left [F,F' \right ] = 
0$;

{\bf (B2)\qua} The map $\pi\circ F$ is  null homotopic in $B\Gg$;
in other words $F=j\circ \al$, where $\al\co S^k\to M$.
\ms

These assumptions permit the construction of
the following commutative diagram:
$$
\CD
 S^k\times S^m @>F\times F'>>        B\Gg_p \\
 @VVV                      @VV\pi V\\
\bigl(S^k\times S^m/S^k\bigr)\;\simeq\;S^{k+m} \vee S^m      @>\{F,F'\}\vee F'>> B\Gg
\endCD
$$
Here $F\times F'$ is some extension of
$F\vee F'$, and exists because $[F,F']=0$.
By (B2) this map descends to
$S^k \times S^m \slash S^k \to B\Gg$.  But
$S^k \times S^m \slash S^k\simeq S^{m+k}\vee S^m$
because the attaching map of the
top cell in the quotient $S^k \times S^m \slash S^k$ is null 
homotopic.  (It is
the Whitehead product of the trivial map on $S^k$ with the identity
map of $S^m$.)
  We denote the
 homotopy class of the induced map on $S^{m+k}$ by
$\{F,F'\}$, 
and by $\{f,f'\}$ the corresponding
 element of $\pi_{k+m-1}(\Gg)$.

Note that the homotopy class $\{F,F'\}$ may depend on the choice of
extension of $F\vee F'$ to $S^k\times S^m$ as well as on the chosen
null homotopy $H_t, t\in [0,1],$ of $\pi\circ F$.  Further it is not
symmetric: indeed  the conditions under which $\{F',F\}$ is defined
are different from those for $\{F,F'\}$.

This
construction is a  particular case of the secondary
Whitehead product.  For any pair of spaces $(X,A)$,
let $\widehat F\in \pi_{k+1}(X,A)$, $F'\in \pi_m(A)$
be such that $[\del \widehat F,F'] = 0 \in \pi_{m+k-1}(A)$. Consider the
following diagram
$$
\CD
\pi_{k+1}(A) @>i>> \pi_{k+1}(X) @>>> \pi_{k+1}(X,A) @>\del>> \pi_k(A)\\
@V[\phantom{F},F']VV @V[\phantom F,F']VV @V[\phantom F,F']VV @V[\phantom F,F']VV\\
\pi_{m+k}(A) @>i>> \pi_{m+k}(X) @>>> \pi_{m+k}(X,A) @>\del>> \pi_{m+k-1}(A)\\
\endCD
$$

Then $[\widehat F,F']\in \pi_{m+k}(X,A)$ has trivial image in
$\pi_{m+k-1}(A)$ and so lifts to $\pi_{m+k}(X)$.  The
secondary Whitehead product
 $$
 \{\widehat F,F'\} \in \pi_{m+k}(X)/i(\pi_{m+k}(A))
 $$
 is defined 
 to be such a lift.
In our situation $X = B\Gg$ and $A = B\Gg_p$.  We are given
two maps $F,F'\in \pi_*(A)$ where $F$ has trivial image in 
$\pi_*(X).$  Therefore we must choose a lift $\widehat F$ of $F$ to
$\pi_{k+1}(X,A)$ which adds an extra indeterminacy $[\pi_{k+1}(X), 
F']$ to the element $\{F,F'\}: = \{\widehat F,F'\}$.

We next specialize further to the case when
$A = BG_p$, where $G_p$ is a Lie group.   Then all Whitehead products 
in $\pi_*(A)\oq$ vanish, and we obtain the following result.

\begin{prop}\label{prop:uniq}  Suppose that the maps $f,f'$
    take values in a Lie group $G_p$ that maps to $\Gg_p$ via $r$.
    Then for each $f\in \pi_{k-1}(G_p)$ such that
    $r\circ f$ is nullhomotopic in $\Gg$
     the above construction gives a homomorphism
     $$
\pi_{m-1}(G_p)\to\pi_{k+m}(B\Gg)\oq/\Kk,\quad
     f'\mapsto \{F,F'\} +\Kk.
     $$
    Here $\Kk$ is the subgroup of $\pi_{k+m}(B\Gg)\oq$ generated by
    $R_*\bigl(\pi_{k+m}(BG_p)\bigr)$
    and the Whitehead products  $[\be,R\circ F']$,  where
    $\be\in \pi_{k+1}(B\Gg)$ and $R\co  BG_p\to B\Gg$ is induced by $r$.
    In particular, if $r\circ f'$ is also nullhomotopic
    in $\Gg$ then $\{F,F'\}$ is  well
    defined modulo the image of $\pi_{k+m}(BG_p)$.  Further in this
    case
    $\{F',F\}$ is defined and equal to $\pm\{F,F'\}$
    modulo the image of $\pi_{k+m}(BG_p)$.
\end{prop}
\begin{proof} The first two statements
    are an immediate consequence of the previous remarks.
    To see that $\{F',F\}=\pm\{F,F'\}$ when $r\circ f$ and $r\circ f'$
    are both null homotopic, observe that
    the two relative Whitehead products 
$$
\begin{array}{ccc}
    (D^{i+j}, S^{i+j-1})&\stackrel W\longrightarrow& (D^{i+1}\vee S^j, S^i\vee 
    S^j)\\
 (D^{i+j}, S^{i+j-1})&\longrightarrow& (S^{i}\vee D^{j+1}, S^i\vee 
    S^j)
    \end{array}
    $$
    become homotopic in  $(D^{i+1}\vee D^{j+1}, S^i\vee 
    S^j)$ since takng the boundary gives an
    isomorphism $\pi_{i+j}(D^{i+1}\vee D^{j+1}, S^i\vee 
    S^j) \to\pi_{i+j-1}(S^i\vee S^j)$.  Hence one gets the same 
    answer if
    one defines the relative product
    $\{F,F'\}$ by using the nullhomotopy of either variable.
    The sign might change when one reverses the order.
\end{proof}

We next give an explicit formula for $\{f,f'\}$ under the
assumptions (A1), (A2).  Choose a map
 $c\co D^{k+m-1}\to \Gg_p$ such that
$\del c = \left<f,f'\right>$ and
an extension $\widetilde f\co D^k\to \Gg$.
Define an element $\be =: \langle\Tilde f, f'\rangle \in 
\pi_{k+m-1}(\Gg)$ as follows.
\begin{equation} \label{eq:ff'}
\begin{array}{cll}
\be: &\bigl(D^k\times S^{m-1}\bigr)\cup D^{k+m-1}\Rar  \Gg,&\\ 
     & \be(x,t) = \left<\widetilde f (x),f'(t)\right> &\text{ for }
      (x,t)\in D^k\times S^{m-1}\\ 
     & \be(z)    =c(z) &\text{ for } z\in D^{k+m-1},
\end{array}
\end{equation}
where $\left<\cdot,\cdot\right>$ is the commutator as before.
In what follows we shall assume that the maps
$f, f'$ take values in $G_p$ and we
choose $c$ also to take values in $G_p$.  Then $\be$
 is well defined modulo the image of $\pi_{k+m-1}(G_p)$
and Samelson products of the form $\left<\ga,f'\right>$, where
$\ga\in \pi_k(\Gg)$.

\begin{prop}\label{prop:repff'}  Let $f, f'$ be as in
    Proposition~\ref{prop:uniq}.
    Denote by $\Kk'$ the subgroup of 
    $\pi_{k+m-1}(\Gg)$ generated by the image of $\pi_{k+m-1}(G_p)$
and Samelson products of the form $\left<\ga,f'\right>$, where
$\ga\in \pi_k(\Gg)$.  Then, if we define $\be = \langle\Tilde f, 
f'\rangle$
by (\ref{eq:ff'}) taking $c$ also to have values in $G_p$, $\be$
represents the  
element of $\{f, f'\}\in 
\pi_{k+m-1}(\Gg)/\Kk'$ corresponding to 
$\{F, F'\}$ of Proposition~\ref{prop:uniq} under the obvious isomorphism $\pi_{k+m-1}(\Gg)/\Kk'\cong \pi_{k+m}(B\Gg)/\Kk$.
\end{prop}

Before giving the proof we derive some easy corollaries.

\begin{cor}\label{cor:rho}  
    Let $\la\co S^1\to \Gg_p$  be a  circle action that is inessential 
in $\Gg$.
    Then $\{\la,\la\}$ is
defined and equals the element
$\rho\in \pi_3(\Gg)$ given by~(\ref{eq:rho3}).  Further, it is 
independent of
choices.
\end{cor}
\begin{proof}  Take $G = S^1$.  Then it is immediate that 
$\{\la,\la\}$
is defined.  It is independent of choices by 
Proposition~\ref{prop:uniq} and
the fact that $\pi_3(G) = 0$.  It equals $\rho$ by 
Proposition~\ref{prop:repff'}.
\end{proof}

\begin{cor} \label{cor:split} Suppose given
a homomorphism $h\co S^1\times K\to \Gg_p$,
where $K$ is a Lie group. Let
$\la: = h|_{S^1\times \{1\}}$ and $f'\in \pi_m(K)$.
Suppose  that $\la$ is inessential in $\Gg$.
Then we may construct the element $\{\la,f'\}\in \pi_{2+m}(\Gg)$
so that it is well
defined modulo the Samelson products
$\left<\be,f'\right>, \be \in\pi_2(\Gg)$.
\end{cor}
\begin{proof}  Define $\{\la,f'\}$ using formula (\ref{eq:ff'})
    noting that we may take $c$ to be constant since, by assumption,
    the maps $\la$ and $f'$ commute.  Then the only indeterminacy is
    of the form $\left<\be,f'\right>,\be \in\pi_2(\Gg)$.
\end{proof}

To show that $\{\la,\la\}\ne 0$ it suffices to show that
$\{\La,\La\}\ne 0$, where $\La\co  S^2\to BS^1$  generates $\pi_2(BS^1)$
and hence that the composite map
$$
S^2\times S^2\Rar BS^1\Rar B\Gg_p \map{\pi} B\Gg
$$
is nontrivial on $H_4$.    This will follow if we show
that the map $BS^1\to B\Gg$ is nontrivial on $H_4$, which we do in
\S\ref{S:ccl}
by evaluating some characteristic classes: see 
Corollary~\ref{cor:nonzero0}.\ms

{\bf Proof of  Proposition~\ref{prop:repff'}}\qua
One can prove this using the following 
commutative diagram that relates the relative Whitehead product to the 
relative Samelson product:
$$
\xymatrix
{
\bigl(D^{i+j-1}, S^{i+j-2}\bigr)\ar[r]^{\Om W\phantom{dupa}} & 
\bigl(\Om(D^{i+1}\vee S^j), \Om(S^{i}\vee S^j)\bigr)\\
\bigl(D^i\times S^{j-1}, S^{i-1}\times S^{j-1}\bigr)\ar[u]\ar[ur]_
{\phantom{dupadu}\langle\Om(\io_1),\Om(\io_2)\rangle}
}
$$
As above, $W$ is the universal model for the relative Whitehead 
product with looping $\Om W$
and the left hand vertical map is induced 
by collapsing the spheres of dimensions $i-1$ and $j-1$. 
Further $\bigl\langle\Om(\io_1),\Om(\io_2)\bigr\rangle$ is the Samelson product 
of the maps obtained by looping the inclusions on each factor.

For readers who are unfamiliar with such homotopy theoretic arguments 
we now give a more explicit proof that describes the relevant 
homotopies in more detail.
We begin by  describing
 the map $\{F,F'\}\co 
S^{k+m}\to B\Gg$.

Consider $S^{k+m}$ as the union of three pieces
$\Dd_1\cup \Dd_0\cup \Dd_{-1}$ where:\ms

{\bf (a)}\qua $\Dd_1$ is the interior of
the top cell in $S^k\times S^m$, and $\{F,F'\}|_{\Dd_1}$ equals
$F\times F'$;

{\bf (b)}\qua $\Dd_0$ is a cylinder $S^{k+m-1}\times [0,1]$ which is
mapped by the
composite
$$
\CD
   S^{k+m-1}\times [0,1]@>e\times id>> \bigl(S^k\vee S^m\bigr)\times
[0,1]
   @>H_t\vee F'>> B\Gg,
\endCD
$$
where $e\co S^{k+m-1}\to S^k\vee S^m$ is the attaching map
of the top cell in $S^k\times S^m$, and $H_t$ is a homotopy from
$H_1 = \pi\circ F$ to the constant map $H_0$; and

{\bf (c)}\qua $\Dd_{-1}$ is a capping disc, mapped by a fixed homotopy
from
$(H_0\vee F')\circ e$ to the constant map. This homotopy exists
because the Whitehead product $[H_0,F']$ vanishes.
Since $F'(S^m)\subset  BG_p$
we may assume that the restriction of $\{F,F'\}$ to $\Dd_{-1}$ factors
through a map $\Dd_{-1}\to BG_p$.  
Thus here we think of the domain as
$$
    S^{k+m}  =  S^{k+m-1}\times
    [-1,2]/\!\sim,
$$
where the equivalence relation $\sim$ collapses each 
sphere 
$$
S^{k+m-1}\times \{-1\},\;\; S^{k+m-1}\times \{2\}
$$
to a single
point and we identify $\Dd_i$ with  $S^{k+m-1}\times [i,i+1]$.

Now choose any model for $B\Gg$ such that the 
induced map $BG_p\to B\Gg$ is injective, and denote by 
$E\Gg$  the space of paths in $B\Gg$ with fixed initial
point.  The map 
 $\Psi\co  S^{k+m}\to \B\Gg$ lifts to
$$
\Tilde\Psi\co D^{k+m}: = D^{k+m-1}\times [-1,2]/\!\sim\;\Rar\; E\Gg.
$$
We must show how to construct $\Psi$  so that 
the restriction of $\Tilde\Psi$ to
the boundary
$$
S^{k+m-1}: = \p D^{k+m-1}\times [-1,2]/\!\sim
$$
agrees with formula (\ref{eq:ff'}).
To do this we construct $\Psi$  using data coming from $\Gg$ rather
than $B\Gg$. Again, we think of  $S^{k+m-1}$ as divided
into three pieces,  a
central cylinder $\Cc_0: = \p D^{k+m-1}\times [0,1]$ together with two
capping discs $\Cc_{-1}, \Cc_{1}$.

Choose the null homotopy $H_t\co  S^k\to B\Gg$, $t\in [0,1]$, to be the 
image under suspension
of the null homotopy $h_t\co S^{k-1}\to \Gg$,
where $h_t$ is the restriction of $\Tf$ to the sphere of radius $t$
in $D^k$.
Then the Whitehead product $[H_t,F']$ desuspends to the Samelson 
product
$\left<h_t,f'\right>$ for each $t\in [0,1]$.  Hence
we may define $\Psi$ on the cylinder $\Dd_0$ to have
a lift $\Tilde\Psi$
that  restricts on $\Cc_0: = \p D^{k+m-1}\times [0,1]$ to the family 
of maps
$\left<h_t,f'\right>$ defined as in~(\ref{eq:sam}).

When $t=0,1$ the Samelson products
$\left<h_t,f'\right>$  lie in $G_p$ and we may choose $\Psi$ and its
lift  $\Tilde\Psi$ so that $\Tilde\Psi$ maps the discs
$D^{k+m-1}\times \{0,1\}$ into  $EG_p$.  Moreover,
$\left<h_0,f'\right>$ is constant, since $h_0$
is the constant map at the identity element.
Thus the restriction of $\Tilde\Psi$ to $D^{k+m-1}\times  \{0\}$
descends to a map $D^{k+m-1}/\p = S^{k+m-1}\to EG_p$.
Then $\Tilde\Psi$ extends to the ball 
  $D^{k+m}$  because $EG_p$ is contractible, and we define
$\Psi$ on $\Dd_{-1}\equiv D^{k+m}$ to be the projection of
$\Tilde\Psi$ to $BG_p$ followed by the map $BG_p\to B\Gg$.

Finally consider the disc $\Dd_1$.
The chosen null homotopy $c$ of
$\left<f,f'\right>$  defines a map of
the disc $D'= \p D^{k+m-1}\times
[1,2]/\!\sim$ to $G_p$ that agrees with the restriction of $\Tilde\Psi$ to 
$\p D' =
\p D^{k+m-1}\times \{1\}$. This gives a map of the $(k+m)$--sphere
$D'\cup \bigl(D^{k+m-1}\times \{1\}\bigr)$ into the contractible
space $EG_p$.  Again, choose $\Tilde\Psi$ on the disc $D^{k+m-1}\times
[1,2]/\!\sim$ to be any extension and define $\Psi|_{\Dd_1}$ to be the
corresponding map to $B\Gg$.
\qed

\subsection{Calculations in $SU(n+1)$}

We end this section by further investigating the construction in 
Corollary~\ref{cor:split} in the case when the homomorphism 
$h$
 takes values in a Lie group $H_p$. Thus we suppose that
 $h\co S^1\times K\to H_p$. 
 For simplicity
we restrict to the case when $H$ is a subgroup of $G: = SU(n+1)$ and 
$H_p$  is contained in 
$G_p:= U(n)$.  (Thus $G/G_p =\cp^n$.)  Unless explicit mention is 
made to the
contrary, we embed $SU(k)$ in $SU(k+1)$ as the subgroup that acts on 
the first
$k$ coordinates, and similarly for inclusions $U(k)\subset U(k+1)$.

For $k = 1,\dots, n$ consider the commuting circle actions
\begin{equation}\label{eq:xx}
\la_k\co S^1\to SU(k+1)\subset SU(n)
\end{equation}
where
$\la_1$ has weights $(1,-1,0\dots,0)$, $\la_2$ has weights 
$(1,1,-2,0,\dots,0)$ and $\la_k$ has weights 
$(1,\dots,1,-k,0,\dots,0)$.
Since $\la_1$ contracts in $SU(2)$ and $\pi_2(U(k))$ $= 0$,
the construction of
Corollary~\ref{cor:split} gives a well defined element
$$
\al_3: = \{\la_1,\la_1\}\in \pi_3(SU(2)).
$$
Because $\la_2$ commutes with $SU(2)$ and contracts in $SU(3)$
we may repeat to get a well defined element
$$
\al_5: = \{\la_2,\al_3\}= \left<\Tilde{\la}_2,\al_3\right> \in 
\pi_5(SU(3)),
$$
where $\Tilde{\la}_2$ is a contraction of $\la_2$ and we 
use formula~(\ref{eq:ff'}).  This is a version of the
construction in Corollary~\ref{cor:split}.  To see this in
general, define the homomorphism
$h_k\co S^1\times U(k)\to U(k)$ to be the identity on the second factor 
and
to take $U(1)$ to the center of $U(k)$ in the obvious way,
 and define $\io_k\co U(k)\to SU(k+1)$ by
$$
\io_k(A) = \left(\begin{array}{cc} A&0\\0&(\det A)^{-1}
                         \end{array}\right).
$$
Then the composite $\io_k\circ h_k\co S^1\to SU(k+1)$ is precisely 
$\la_{k}$,
and assuming that $\al_{2k-1}\in \pi_{2k-1}(SU(k))$ is already defined
we may set:
$$
\al_{2k+1}: = \{\la_k,\al_{2k-1}\} = 
\left<\Tilde{\la}_k,\al_{2k_1}\right>
\in \pi_{2k+1}(SU(k+1)),
$$
where $\Tilde{\la}_k\co D^2\to SU(k+1)$ is a contraction of $\la_k$.

\begin{lemma}\label{le:un}  
$\al_{2k+1}\in\pi_{2k+1}(SU(k+1))$ is nonzero for $k= 1,\dots n$.
\end{lemma}
\begin{proof} 
    Let $\T^n$ be the diagonal torus in $SU(n+1)$ and
denote by $\La_k\co S^2\to B\T^k$ the desuspension of $\la_k$.
Because Whitehead products vanish in $B\T^k$, the map
 $\La_1\vee\La_1\vee\La_2\vee\dots\vee\La_k\co \vee_{k+1} S^2\to B\T^k$ 
 has an extension to 
 $$
 \CD
 \prod_{i=1}^{k+1} S^2 
@>g_k=\La_1\times\La_1\times\La_2\times\dots\times\La_k>>
 B\T^k
 \endCD
$$
that is unique up to homotopy.  Taking the composite with
 the inclusion $B\T^k\to 
BU(k+1)$ gives a rank $(k+1)$ bundle over $\prod S^2$. We show below 
that it suffices to show that this bundle has some nontrivial Chern 
classes in dimension $2k+2$.  We then calculate the pullback of $c_{k+1}$.

We first claim that for each $k\in \{1,\dots,n\}$ 
there is a homotopy commutative diagram
$$
\CD
\prod_{i=1}^{k+1} S^2 
@>g_k>>
 B{\T}^k\\
@VV\ga_{k} V @VVV\\
S^{2k+2}\vee J_k @>F_{k}>> BSU(k+1),
\endCD
$$
in which 
the CW complex $J_k$ 
has dimension $\le 2k$.
We construct this diagram by induction on $k$.
The only difficulty is to construct the left vertical map $\ga_k$.
When $k=1$ the diagram reduces to
$$
\CD
S^2\times S^2 
@> \La_1\times \La_1>>
 B\T^1\\
@VV\ga_1 V @VVV\\
 S^4\vee S^2@>F_1>> BSU(2)\subset \rlap{$BSU(n+1),$}
\endCD
$$
and is a special case of the situation discussed before
Proposition~\ref{prop:uniq}.  For $k>1$
we may assume by induction
that the map $\prod_{i=1}^{k+1} S^2\to BSU(k+1)$ given by going 
horizontally
to $B\T^k$ and then vertically to $BSU(k+1)$  factors through
$$
\CD
S^2\times \bigl(S^{2k}\vee J_k\bigr) @>\La_k\times F_{k-1}>> 
BZ_k\times BSU(k)\to BSU(k+1) 
\subset BSU(n+1),
\endCD
$$
where  $Z_k$ denotes the commutator of $SU(k)$ in $SU(k+1)$.
Since $\La_k$ contracts in $BSU(k+1)$ the product
 $\La_k\times F_{k-1}$ is homotopic via maps to $BSU(k+1)$
to a map that takes the sphere $S^2\times \{pt\}$ 
to a single point.  But when one contracts the first sphere
in the product $S^2\times (S^{2k}\vee J_{k-1})$ 
one obtains a space of the form $S^{2k+2}\vee J_k$.
Thus the above diagram exists for all $k$.

Next observe that  $g_k$
factors through the inclusion $BU(k)\to BSU(k+1)$
given by
$$
A\mapsto \left(
\begin{array}{cc}
A & 0\\
0 & \det A^{-1}
\end{array}
\right).
$$
Therefore we may apply 
Proposition~\ref{prop:repff'} to conclude that
 the restriction of $F_k$ to the top sphere $S^{2k+2}$ 
desuspends to $\al_{2k+1}\in \pi_{2k+1}SU(k+1))$ modulo the 
subgroup
$\Kk'$ generated by the image of 
 $\pi_{2k+1}(U(k))$ and certain Samelson products.  Since $\Kk' = 0$,
we find that $\al_{2k+1}\ne 0$ if and only if
$F_k|_{S^{2k+2}}$ is nonzero.  
This will be the case precisely when the bundle represented 
by
the composite
$$
\CD
\prod S^2@>g_k>>
B\T^{k} @>>>BSU(k+1)
\endCD
$$
has nontrivial top dimensional Chern classes. 

We now check this by calculating the pullback of $c_{k+1}$.
\footnote
 {
Because the map $BSU(k+1)\to 
 BU(\infty)$ is $(2k+2)$--connected, one can equally well
phrase this calculation in terms of the structure of the induced 
stable bundle on $\prod S^2$.} 
Denote by $\T: = (S^1)^{k+1}$ the diagonal subgroup in $U(k+1)$
and by $t_1,\dots, t_{k+1}\in H^2(B\T)$  the obvious generators of $H^*(B\T)$.
Then $c_{k+1}$ pulls back to $t_1\dots t_{k+1}\in H^{2k+2}(B\T)$.
The map $B\T^{k} \to BU(k+1)$ factors through $B\T$ and it suffices
to show that the pull back of $t_1\dots t_{k+1}$ by the composite
$$
\CD
\Phi\co  \prod S^2@>g_k>>
B\T^{k} @>>>B\T
\endCD
$$
does not vanish.  
For $j = 0,1,\dots,k$ denote by
 $y_j$  the pullback to $\prod_{j=0}^{k} S^2$ of the 
 generator of $H^2(S^2)$ by the projection onto 
 the sphere that is mapped by $\La_j$
 (where we set $\La_0: = \La_1$.)  Then 
 \begin{eqnarray*}
 \Phi^*(t_1\dots t_{k+1}) &=& \Bigl(y_0+y_1+y_2 +\dots + y_k\Bigr) \Bigl(
 -y_0-y_1+y_2 +\cdots + y_k\Bigr)\times \\
&&\qquad\quad \times 
\prod_{j=2}^k\Bigl(-jy_j + \sum_{i>j}y_i\Bigr)\\
& = & 2(-1)^kk!\,y_0y_1y_2\dots y_k.
 \end{eqnarray*}
since $y_j^2= 0$ for all $j$.
 Since this is nonzero, the proof is complete. 
\end{proof}

\begin{remark}\rm
One can construct a  nontrivial element in
$\pi_{2k+2}(BSU(k+2))$ for $k\ge 1$ by the following inductive 
procedure.
 Let $E_2\to S^2$ be the complex
line bundle classified by the map $\La_2\co S^2\to BU(1)$. Note that
its first Chern class is equal to the generator $y\in H^2(S^2)$.
Suppose also that we have already constructed
a map $\La_{2k}\co S^{2k}\to BU(k)$ that classifies a bundle $E_{2k}$
with nontrivial Chern class $c_k(E_{2k})$.
Then construct a homotopy commutative  diagram
$$
\xymatrix{
S^2\times S^{2k} \ar[d] \ar[rr]^{\!\!\!\!\!\!\La_2\times \La_{2k}}& & 
BU(1)\times BU(k)\ar[r]
&BU(k+1)\ar[d]\\
S^{2k+2}\vee S^{2k}\ar[rrr]^{F} & & & BSU(k+2)
}
$$
in which the top row classifies  the bundle
$E_2\times E_{2k}\to S^2\times S^{2k}$ whose Chern class
$c_{k+1}=p_1^*c_1(E_2)\cup p_2^*c_k(E_{2k})\neq 0$.
Since $\La_2$ contracts in $BSU(k+1)$ the bottom row of the diagram 
is constructed in the usual way, and the induced map $S^{2k+2}\to 
BSU(k+2)$
is homotopically nontrivial.  But note that it does {\it not} 
desuspend to
the element of $\pi_{2k+1}(SU(k+2))$ constructed via commutators
as in Corollary~\ref{cor:split}.  For one can choose
the contraction of $\La_2$ to lie in a copy of $SU(2)\subset SU(k+2)$ 
that {\it commutes} with the image of $U(k)$.  
This is not a contradiction since Proposition~\ref{prop:repff'}
states only that the induced map should be given by~(\ref{eq:ff'})
modulo elements of $\Kk'$.  But now $\Kk'$ is the whole of
$\pi_{2k+1}(SU(k+2))$ since it includes the image of 
$\pi_{2k+1}(U(k+1))$.
\end{remark}

\section{Characteristic classes}\label{S:ccl}

Let $\Gg=\Ham\Mo$  and consider the universal bundle
$M\stackrel{j}\to M_\Gg \stackrel{\pi}\to B\Gg$.
There is a unique class $[\Om]\in H^2(M_\Gg;\R)$
called the {\em coupling class} that extends the fiberwise symplectic
class
$[\om]$ and has the property that the fiberwise integral
$\int_M[\Om]^{n+1}\in H^2(B\Gg;\R)$ vanishes.   Following
Januszkiewicz--K\c edra~\cite{JK}, we define the classes
\begin{equation}\label{eq:ka0}
\mu_k: = \int_M [\Om]^{n+k}\in H^{2k}(B\Gg): = H^{2k}(B\Gg;\R).
\end{equation}
In this section we first 
generalize  these classes to other groups and monoids, 
and then discuss computations and applications.  
When $H^1(M)\ne 0$ we shall work only with
 {\it connected} groups and 
monoids.  The issues that arise in the general case are discussed
from different perspectives in
Gal--Kedra~\cite{GalK} and McDuff~\cite{Me}.

\subsection{The classes $\mu_k$}\label{ss:muk}

Our first aim is to define the classes $\mu_k$ in as general a context 
as possible.  Thus if $\Hh$ is a topological monoid that acts on a
c-symplectic manifold $\Ma$ we need to determine conditions on $\Hh$ 
that guarantee that the class $a\in H^2(M)$ has a well defined 
extension to a class $\Tilde a\in H^2(M_{\Hh})$.
As before we denote by $\Hh_0$  the identity component of the space of 
smooth maps $M\to M$
and by 
$
\Hh_a
$
the space of smooth homotopy equivalences that fix 
$a$.  Further $M_{\Hh}$ denotes the total space of the universal 
$M$--bundle over $B\Hh$.  Thus the projection $\pi\co  M_\Hh\to B\Hh$
is a Hurewicz fibration.

The $a$--Flux homomorphism $\Flux^a\co  \pi_1(\Hh_0)\to H^1(M;\R)$
was defined in equation~(\ref{eq:aflux}).
In the symplectic case the image 
$$
\Ga_\om: = \Flux^a\bigl(\pi_1(\Symp_0)\bigr)\subset H^1(M;\R)
$$
is called the flux group, and
there is  a surjective homomorphism 
\begin{equation}\label{eq:flux}
\Flux_\om\co\Symp_0\Mo\to H^1(M;\R)/\Ga_\om,\qquad
g\mapsto \int_0^1[\om(\dot g_t,\cdot)]dt,
\end{equation}
where $g_t$ is any path in the connected group
$\Symp_0\Mo$ from the identity to $g$.  The  kernel 
of $\Flux_\om$  
is precisely the Hamiltonian group.
In the c-symplectic case one cannot define a flux homomorphism on
$\Diff_0(M)$ (which is a simple group),
and the analog 
of the Hamiltonian group is the covering group
$$
\Ham\Ma
$$
of $\Diff_0(M)$ corresponding to the kernel of
 $
 \Flux^a\co  \pi_1(\Diff_0)\to H^1(M;\R).
 $
 Thus there is an exact sequence of topological monoids
$$
 \Ga_a \to \Ham\Ma\to \Diff_0M,
 $$
 where the fiber $  \Ga_a:= {\rm im}(\Flux^a)$ is given the discrete topology.
 \footnote
{Ono~\cite{Ono} recently proved that the symplectic flux group $\Ga_\om$ is a discrete subroup of $H^1(M;\R)$.  But $\Ga_a$ need not be: see McDuff~\cite{Mcf}.}
Note that in the symplectic case, $\Ham\Mo$ is 
 homotopy equivalent to  the corresponding covering group of $\Symp_0\Mo$. 
 Further there is an inclusion $\Ham\Mo\to \Ham\Ma$ given by taking 
 the element $f\in \Ham\Mo$ to the pair $[f, \{f_t\}]$, where $f_t, 
 t\in [0,1],$ is any path in the connected group $\Ham\Mo$ from the 
 identity to $f$.  (Since two such paths differ by a loop with zero 
 $a$--flux, this recipe defines a unique element $[f, \{f_t\}] \in 
 \Ham\Ma$.) 
Sometimes we will also work with
the covering space 
$
\THha
$
of  $\Hh_0$ 
corresponding to the kernel of 
$
 \Flux^a.
$
This is a topological  monoid.

A smooth circle action on a c-symplectic manifold $(M,a)$
is said to be {\em c-Hamiltonian} if it is in the kernel of
$\Flux^a$.  As shown by Allday~\cite{All}
such actions do have certain geometric properties; for example their
fixed point set has at least two connected components, though these
components need not  be c-symplectic  as would happen in the
symplectic case.

\begin{prop}\label{prop:ka} {\rm (i)}\qua Given $a\in H^2(M)$ such that 
$\int_Ma^n\ne 0$,
consider the fibration $M\to M_{\Hh}\to B\Hh$ where $\Hh: = \THha$.
Then there is a unique
element $\Tilde a\in H^2(M_{\Hh})$ 
that restricts to $a\in H^2(M)$
and is such that $\int_M \Tilde{a}^{n+1} = 0$.

{\rm (ii)}\qua   If $H^1(M)=0$,
such a class $\Tilde a$  exists on
$M_{\Hh_a}$.
\end{prop}

The above class $\Tilde a$ is called the {\em coupling class}.

\begin{proof}(i)\qua  
Consider the Leray--Serre cohomology spectral sequence
 for the (Hurewicz)
fibration $M\to M_{\Hh}\to B\Hh$.  Since $\Hh$ is connected the
$E_2$--term is a product and the class $a$ lies in $E_2^{2,0} = H^2(M)\otimes
H^0(B\Hh)$.  As explained
for example in~\cite[Lemma~2.2]{LM} the differential
$$
d_2\co  E_2^{2,0} \Rar E_2^{1,2} =
H^1(M)\otimes H^2(B\Hh)
$$
is determined by the flux homomorphism.  More precisely,
if we consider the elements of $H^1(M)\otimes H^2(B\Hh)$ as
homomorphisms $H_2(B\Hh)\to H^1(M)$, then
$$
d_2(a)(\ov{\la}) = \Flux^a(\la),
$$
where $\ov\la\in \pi_2(B\Hh)\cong H_2(B\Hh)$ corresponds to the loop
$\la\in \pi_1(\Hh)$.  Hence, because $\la$ is c-Hamiltonian,
$d_2(a) = 0$.
But
the image of $a$ under $d_3\co  E_3^{2,0} \to E_3^{0,3}$ must vanish:
since $a^{n+1} = 0$ in $H^*(M)$,
$$
0 = d_3(a^{n+1}) = (n+1)a^n \otimes d_3(a) \in H^{2n}(M)\otimes 
H^{3}(B\Hh),
$$
which is possible only if $d_3(a) = 0$.
Therefore $a$ survives into the $E_{\infty}$--term of the spectral
sequence, and so has some extension $u\in H^2(M_{\Hh})$.

To prove  uniqueness note that because $B\Hh$ is simply connected the
kernel of the restriction map
$H^2(M_{\Hh})\to H^2(M)$ is isomorphic to the  pullback of 
$H^2(B\Hh)$.
A short calculation shows that we may take
$$
\Tilde a : = u -\frac 1{n+1}\pi^*\pi_{!}\bigl(u^{n+1}\bigr),
$$
where $\pi_!$ denotes integration over the fiber.
This proves (i).

Now consider (ii).
 Because the elements in $\Hh: = \Hh_a$ preserve $a$
there is an element $a$ in the $E_2^{0,2}$--term in the Leray--Serre 
spectral sequence for the Hurewicz fibration
$M_{\Hh_a}\to B\Hh_a$.  Since $H^1(M) = 0$ the differential
 $d_2$ must vanish and the argument showing that
$d_3(a)=0$ still holds.  
Further the uniqueness 
proof  goes through as before: although $B\Hh$ need no longer be simply 
connected, the fact that $H^1(M)=0$ implies that $E_*^{1,1} = 0$.  
Hence the
kernel of the restriction map
$H^2(M_{\Hh})\to H^2(M)$ is  still isomorphic to the  pullback of 
$H^2(B\Hh)$.  This proves (ii).
\end{proof}

\begin{rmk}\rm (i)\qua If the  first Chern class of $\Mo$ is a 
nonzero multiple of the symplectic class $a$, then 
there is an easier way to find an
extension $u\in H^2(M_{\Hh})$ of the symplectic class: simply take 
it to be an appropriate multiple 
of the first Chern class of the vertical tangent bundle. This 
construction applies  whenever
$\Hh$ acts in such a way that the vertical bundle has a complex structure.
    The most natural choice for $\Hh$ is the group of 
    symplectomorphisms.  Note that we cannot take it to be the diffeomorphism 
    group.\medskip
    
  (ii)\qua  Let $(M,a)$ be any c-symplectic manifold and consider an 
 $M$--bundle $P\to B$.
 We saw above that if $a$ survives into the $E_2$ term of 
 the Leray--Serre spectral sequence then 
 $d_3(a) = 0$ and $a$ survives to $E_{\infty}$.  This shows that the 
 obstruction to the existence of an extension $\Ta$ of $a$ depends 
 only on the restriction of $P\to B$ over the $2$--skeleton of $B$. 
    \end{rmk}

\begin{defn}\label{def:ka}\rm Let $(M,a)$ be a c-symplectic manifold and let
    $\Hh$ denote either the monoid 
$\THha$  or, if $H^1(M)=0$, the monoid $\Hh_a$.
We define $\mu_k\in H^{2k}(B\Hh)$ by 
$$
\mu_k:= \pi_!(\Tilde {a}^{n+k}),
$$
where $\Tilde a$ is the coupling class 
constructed in Proposition~\ref{prop:ka}.
\end{defn}

These classes $\mu_k$ extend those defined in
equation~(\ref{eq:ka0}).  Note that $\mu_1\equiv 0$ by definition.
We will see in the next section that the pullback of 
$\mu_{2k}$ to $H^*(BS^1)$ is
nonzero for all $k>0$ and for every 
nontrivial c-Hamiltonian circle action on $M$.\ms

{\bf Proof of Proposition~\ref{prop:cpn}}\qua
    We must show that if $G = SU(n+1)$ acts on $M:=\cp^n$
    then the induced map $R_*\co \pi_*(BSU(n+1)\oq)\to H_*(B\Hh)$  is
injective where 
$\Hh: = \Hh_a$.
For $2\leq k \leq n+1$  choose $s_{2k} \co S^{2k}\to BSU(n+1)$
so that
$\al_k: = s_{2k}^*(c_k)$ generates
$H^{2k}(S^{2k})$,
where $c_k$ is the $k$th Chern class.
Consider the bundle $\cp^n\to P\stackrel{p}\to S^{2k}$
associated with this element.  (This is simply the projectivization 
of the corresponding rank $(n+1)$ vector 
bundle $E\to S^{2k}$.)  The Leray--Hirsch theorem
states that the cohomology of $P$ is a free $H^*(S^{2k})$--module 
generated by the powers $1=c^0,\dots,c^n$ of the first Chern class $c$
of the tautological line bundle.  
Since $2\le k\le n+1$ there is $\be_k\in H^{2k}(S^{2k})$ so that
$$
c^{n+1} = p^*(\be_k)\cup c^{n+1-k}.
$$
The theory of 
characteristic classes implies that $\be_k = c_k(E)= s_{2k}^*(c_k)=:\al_k.$
Multiplying this equality by $c^{k-1}$ we obtain
$$
c^{n+k} = p^*(\al_k)\cup c^n,
$$
which implies that $p_!(c^{n+k})= const.\al_k\ne 0$.  Now observe that 
because $k\ge 2$, $H^2(P)$ has 
dimension $1$, so that 
$c$  equals  the coupling class $\Tilde a\in H^2(P)$ up to 
a nonzero constant.
Hence 
$$
\left <\mu_k, R_*([S^{2k}])\right > =
const. \left<\al_k,[S^{2k}]\right >\neq 0.
$$
Therefore $\pi_*(BG)\to H_*(B\Hh)$  is injective. The conclusion now follows as in the proof of Corollary~\ref{cor:cfl}.
\qed

\begin{rmk}\label{rmk:p}\rm
Let $\Gg: = \Ham\Ma$ and denote by $\Gg_p$ the homotopy fiber of the 
evaluation map $\Gg\to M$.  Thus if the elements of $\Gg$ are pairs 
$(g, \{g_t\})$ where $\{g_t\}$ is an equivalence class of paths from 
$id$ to $g$, then $\Gg_p$ is the subgroup of $\Gg$
consisting of all pairs such 
that $g(p) = p$.
The restriction $M_{\Gg_p}\to B\Gg_p$
of the universal 
bundle $M_{\Gg}\to B\Gg$ to $B\Gg_p\subset B\Gg$ has a canonical 
section $\si\co  B\Gg_p\to M_{\Gg_p}$ whose image $\si(b)$ at $b\in 
B\Gg_p$ is the  point 
in the fiber $M_b$ corresponding to $p$.  
Thus there is a homotopy commutative diagram
$$
\xymatrix
{
M_{\Gg_p}\ar[r]^\io \ar[d]& M_\Gg\ar[d]\\
B\Gg_p \ar@/^1.1pc/[u]^{\si}\ar[r] & B\Gg
}
$$
such that the composite $\io\circ \si\co  B\Gg_p\to M_{\Gg}$ is a 
homotopy equivalence.
By Proposition~\ref{prop:ka}, the fiberwise symplectic class $a$
extends to $M_{\Gg_p} \subset M_{\Gg}$.
We shall denote by $\Tilde a_p$ the 
extension that is normalized by the requirement that $\si^*(\Tilde 
a_p) = 0$.  Correspondingly there are characteristic classes
$$
\nu_k: = \pi_!\Bigl((\Tilde a_p)^{k+n}\Bigr)\in H^{2k}(B\Ham\Ma_p),\quad k\ge 1.
$$
We show at the end of \S\ref{ss:circ} that $\nu_1$ need not vanish. 
\end{rmk}

\subsection{Calculations for circle actions}\label{ss:circ}

We suppose that $S^1$ acts smoothly on a connected  almost symplectic
manifold $(M,a)$.  By averaging we may construct an $S^1$--invariant
closed representative $\om$ of the class $a$.  If $\underline \xi$ is 
the
generating vector field on $M$ for the action then
the identity
$$
0 =\Ll_{\underline \xi}(\om) = d(\io_{\underline \xi}\om) + 
\io_{\underline
\xi}(d\om)
$$
implies that the $1$--form $\io_{\underline \xi}\om$ is closed.  To say
the action is c-Hamiltonian is equivalent to saying that this
$1$--form is exact.  The requirements
$$
dH = \io_{\underline \xi}\om,\quad  \int_M H\om^n =0
$$
define a unique function $H\co M\to \R$ that is called the normalized
(Hamiltonian) $\om$--moment map.  The following useful result is 
presumably
well known.

\begin{lemma}\label{le:cham}
If the action is nontrivial (ie nonconstant) $H$ cannot be 
identically
zero.
\end{lemma}

\begin{proof}\, There is a subset $M_0$ of full measure
in $M$ where the $S^1$ action gives rise to a fibration $S^1\to
M_0\to N$.  If $\om$ is both invariant and such that
$\io_{\underline \xi}\om \equiv 0$ then $\om|_{M_0}$ pulls back from 
$N$ and
so cannot satisfy the condition $\int_M\om^n \ne 0$.
\end{proof}

We first describe a cocycle representing the coupling class
in the Cartan model of the equivariant cohomology
(see \cite{gs} for details).
Recall that the Cartan model is the following DGA
$$
\Om^*_{S^1}(M) := \left (S(\Lie(S^1)^*)\otimes 
\Om^*(M)\right)\,\!^{S^1},
$$
with differential
$d:= 1\otimes d_M  - x \otimes \iota_{\underline \xi}$,
where  $\xi \in \Lie(S^1)$ is a basis vector, $x\in S(\Lie(S^1)^*)$ 
its dual
and $\underline \xi$ is the induced vector field on $M$.
The coupling class is represented by
$\Om :=1\otimes \om - x \otimes H$, where $H\co M\to \B R$ is
the normalized (ie $\int_M H\om^n =0$) $\om$--moment map.
Note that, by definition of the moment map,
$$
d\Om = -x \otimes \iota_{\underline \xi}\om - x \otimes d_MH = 0.
$$
Thus $\Om$ is closed.

The fiber integration
$\pi_!\co \Om^*_{S^1}(M)\to \Om^*_{S^1}(pt) = S(\Lie(S^1)^*)$
corresponds to the  (equivariant) constant  map $M\to pt$ in
the following way. An equivariant differential form
$\al \in \Om^*_{S^1}(M)$
can be regarded as a polynomial map $\al \co \Lie(S^1)\to \Om^*(M)$.
Then $\pi_!(\al)$ is a polynomial map $\Lie(S^1)\to \B R$ given
by
$$
(\pi_!(\al))(\xi ) := \int_M \al(\xi),
$$
since the (nonequivariant) fiber integration corresponding
to the constant map is just the usual integration over the manifold
$\int_M\co \Om^*(M)\to \Om^*(pt)=\B R$ (cf. Theorem 10.1.1 in \cite{gs}).

\begin{lemma}\label{le:fint}  Suppose that
$(M,a)$ is a c-symplectic manifold 
of dimension $2n$ and consider a c-Hamiltonian circle action $\la$
on $M$ with  normalized 
$\om$--moment map $H$.   Denote by $\mu_k(\la)\in H^*_{S^1}(pt)$ the 
pullback of the 
characteristic class $\mu_k$ to $BS^1$ by the classifying map of the
 associated bundle $M\to P\to BS^1$.  
 Then
\begin{equation}\label{eq:ck}
    \mu_k(\la)= (-1)^k
\left(\begin{array}{c}n+k \\ n \end{array}\right)
\int_M H^k\om^n \cdot x^k \in S(\Lie(S^1)^*)=H_{S^1}^*(pt).
\end{equation}
In particular, 
$\mu_k(\la)\ne 0$ for even $k$ whenever the
circle action is nontrivial.
\end{lemma}

\begin{proof}
Because $\om^m = 0$ for $m>n$ for dimensional reasons,
\begin{eqnarray*}
\Om^{n+k} &=& (1\otimes \om - x\otimes H)^{n+k} \cr
         &=& \sum _{m=0}^{n} (-1)^{n+k-m}
\left(\begin{array}{c} n+k \\ m \end{array} \right)
\left(x^{n+k-m}\otimes \om^m H^{n+k-m}\right).
\end{eqnarray*}
Now apply the above discussion to see that the only term
that contributes to the integral is the one containing the $n$-th power
of the symplectic form.  This establishes formula~(\ref{eq:ck}).
Now observe that the function $H^k$ is nonnegative for even $k$ and
so by Lemma~\ref{le:cham} its integral vanishes only if the action is trivial.
\end{proof}

We now complete the proof of
Theorem~\ref{thm:s0}.  

\begin{cor}\label{cor:nonzero0}
    Let $\Hh$ be as in Definition~\ref{def:ka}
 and suppose that  $\la\subset \Hh$ is an inessential
    (nontrivial) circle
    action on $M$.  Denote by 
        $\ov{\rho}\in \pi_4(B\Hh)$ the element $\{\La,\La\}$ that is 
        formed from the suspension $\La\co S^2\to B\Hh$ of $\la$ as in 
        Proposition~\ref{prop:uniq}.
Then
    $$
        \mu_2(\ov{\rho})\ne 0.
    $$
\end{cor}

\begin{proof} Denote by $r\co S^1\to \Hh$ the homomorphism with
image $\la$.
As remarked after Corollary~\ref{cor:split} it suffices
to show that the corresponding map $R\co  BS^1\to B\Hh$ is nontrivial
on $H_4$.  But it follows from Lemma~\ref{le:fint} that the pullback
$R^*(\mu_2) = [\pi_!(\Om^{n+2})]$ is nonzero, where $\mu_2$ is as
in~(\ref{eq:ka0}).
\end{proof}

The next result is an immediate consequence of Lemma~\ref{le:fint} and 
the proof of Lemma~\ref{le:un}.

\begin{lemma}\label{le:nonzero1}
    Suppose that $G: = SU(\ell)$ acts smoothly on 
    the $2n$--manifold $(M,a)$, and  let $\Hh$ be as in Definition~\ref{def:ka}.
    Denote by $\la_k, k = 1,\dots,\ell-1$, the circles in 
    $SU(\ell)$ defined by~(\ref{eq:xx}) and by $H_k, k = 
    1,\dots,\ell-1$,
    the $\om$--moment maps for $\la_k$ where $\om$ is a $G$--invariant 
    representative for $a$.  Then 
    $
    \pi_{2k}(BG)$ has nonzero image in $H_{2k}(B\Hh)$
    if and only if 
    $$
    \int_M (H_1)^2H_2\cdots H_{k-1} \om^n\ne 0.
    $$
\end{lemma}

{\bf Proof of Proposition~\ref{prop:nonzero}}\qua
This is the special 
     case with $\ell = 2$.  Then the integral is $\int H_1^2\om^n$ 
     which is always nonzero by Lemma~\ref{le:cham}.
\qed\ms

{\bf Proof of Lemma~\ref{le:p}}\qua  Consider a c-Hamiltonian $S^1$--action 
with moment map $H$, and choose a fixed point $p$.
Because  $S^1$ maps to a loop in $\Diff_0$ with trivial flux
the inclusion $S^1\to \Diff_0$ lifts to a homomorphism $S^1\to\Gg_p$,
where $\Gg_p$ is as defined in 
Remark~\ref{rmk:p}.  Moreover, the section $\si$ pulls back over 
the corresponding classifying map $BS^1\to B\Gg_p$ to the section
of $ M_{S^1}\to BS^1$ with image $ BS^1\times\{p\}$.  Therefore the
pullback of $\Tilde a_p$ to $M_{S^1}$ may be 
represented in the Cartan model by the element $\Om_a: = 1\otimes 
\om - x\otimes (H-H(p))$.  Hence the proof of Lemma~\ref{le:fint} 
shows that for some  nonzero constant $c$
$$
\nu_1(\la) = c. \int_M(H-H(p))\om^n.
$$
If this integral is nonzero, the $2$--skeleton $S^2$ of $BS^1$ maps 
nontrivially to  $B\Gg_p$ under the
classifying map $BS^1\to B\Gg_p$.  This is equivalent to 
saying that $\la$ is essential in $\Gg_p$.\qed

\subsection{Relation to Reznikov's classes and the 
$\ka$--classes}\label{ss:rez}

Above we defined classes $\mu_k$ in $H^*(B\Hh)$ for monoids $\Hh$
whose action is c-Hamiltonian, ie, the c-symplectic class $a$
extends over $M_{\Hh}$.  If we restrict to the case $\Gg = \Ham\Mo$
then these classes have a geometric
interpretation.  We now show that they can be constructed by the usual
Chern--Weil process.   As a consequence they desuspend
to classes in $H^{2k-1}(\Ham\Mo)$ 
(the usual singular cohomology of $\Ham\Mo$ considered as a topological space) that are the restrictions of 
Reznikov's
classes in the smooth Lie algebra cohomology 
$H_{sm}^{2k-1}(\Lie\Ham\Mo)$.

To understand this, recall that  $\Lie\Ham\Mo$ 
 can be identified with the space of functions
$C_0(M)$ of zero mean on $M$ (with respect to the volume form
$\om^n$.)   Moreover, it has a nondegenerate bilinear form (a Killing
form) defined by
$$
\langle F,G\rangle : = \int_M FG \;\om^n
$$
that is invariant under the adjoint action of $\Gg: = \Ham\Mo$.  Now
the coupling class $[\Om]$ may be represented
by a closed differential form $\Om$  on $M_\Gg$ that  extends the 
fiberwise
symplectic form.  (Here we are working with de Rham theory on $M_\Gg$
which is not, strictly speaking, a manifold.  However, one can make
everything precise by considering an arbitrary smooth map $B\to B\Gg$
and looking at the associated smooth pullback bundle over $B$.)
The form
$\Om$ defines a connection on the bundle
$\pi\co  M_\Gg\to B_\Gg$ whose horizontal distribution $Hor_x, x\in 
M_{\Gg},$
is given by the  $\Om$--orthogonals to the vertical tangent bundle.
As shown by Guillemin--Lerman--Sternberg~\cite{GLS}, the corresponding
parallel transport maps  preserve the fiberwise symplectic form and 
the
holonomy is Hamiltonian.  Moreover,
given vector fields $v,w\in T_bB\Gg$ with horizontal lifts
$v^{\sharp}, w^{\sharp}$, the function $\Om(v^{\sharp}, 
w^{\sharp})(x)$
restricts on each fiber $M_b: = \pi^{-1}(b)$ to an element of
$\Lie(\Gg)$ that represents the curvature of this connection at 
$(v,w)$.
In other words, the closed $2$--form
$\Om$ on $M_\Gg$  defines a $2$--form $\Tilde\Om$
on the base $B\Gg$ of the
fibration that takes values in the Lie algebra $\Lie(\Gg)$
of the structural group $\Gg$ and is a curvature form
in the usual sense.  (For more detail, see
McDuff--Salamon~\cite[Ch~6]{MS}.)

Any Ad-invariant polynomial
$\Ii^k\co  \Lie(\Gg)^{\otimes k}\to
\R$ therefore gives rise to a
characteristic class $c_k^\Gg$ in $H^*(B\Gg)$, namely the class 
represented by
the closed real-valued $2k$--form $\Ii^k\circ \Tilde{\Om}^k$.  In the 
case at
hand we may take
$$
\Ii^k (F_1\otimes\dots\otimes F_k): =\int_M F_1\cdots F_k
\;\om^n.
$$

\begin{lemma}  This class $c_k^\Gg$ equals $const. \mu_k$.
\end{lemma}
\begin{proof}\,
Let $v_1,\dots,v_{2k}$ be vector fields on $B\Gg$ with
horizontal lifts $v_1^{\sharp},\dots,v_{2k}^{\sharp}$.  Then,
if the $w_j$ are tangent to
the fiber at $x\in M_{\Gg}$ we find
\begin{eqnarray*}
\Om^{n+k}(w_1,\dots,w_{2n},v_1^{\sharp},\dots,v_{2k}^{\sharp})(x)
& = & \sum_\si \eps(\si) \left(\begin{array}{c}n+k \\ n 
\end{array}\right)\;\;\times \\
&& \;\;  F_{1,\si}(x)\cdots F_{k,\si}(x)\;
\om^n(w_1,\dots,w_{2n}),
\end{eqnarray*}
where, for each permutation $\si$ of $\{1,\dots,2k\}$, $\eps(\si)$
denotes its signature and
$$
F_{j,\si}(x): =
\Om(v_{\si(2j-1)}^{\sharp},v_{\si(2j)}^{\sharp})(x) =
\Tilde\Om(v_{\si(2j-1)},v_{\si(2j)})(x).
$$
Therefore
$\bigl(\pi_!\Om^{n+k}\bigr)(v_1,\dots,v_{2k}) = const.
\Ii^k\circ \Tilde{\Om}^k(v_1,\dots,v_{2k})$ as claimed.
\end{proof}

Now suppose given a homomorphism from a Lie group $G$ to $\Ham\Mo$.
The classes $\mu_k\in H^*(B\Ham\Mo)$ pull back under the map $R\co BG\to
B\Ham\Mo$ to {\it some} elements in $H^*(BG)$, which by definition
are $G$--charact\-eristic classes.  One can figure out which classes
they are by investigating the invariant polynomial
$$
\Lie(G)^{\otimes k}\Rar C_0(M)^{\otimes k}\stackrel{\Ii^k}\Rar \R.
$$
Reznikov~\cite{Rez} did this calculation (in a slightly different
context) for the case of the action of $SU(n+1)$ on $\cp^n$, and 
concluded that the $\mu_k$ pull back to algebraically independent 
elements  in $H^*(BSU(n+1))$.
The advantage of our approach is that the classes extend to 
$B\THha$ (and when $H^1(M)=0$ to $B\Hh_a$)
since one does not use the action of the Lie algebra in their 
definition.

There are other characteristic classes in $H^*(B\Gg)$  arising from 
the characteristic classes on $M$  preserved by the elements 
of $\Gg$, ie the Chern classes when $\Gg= \Symp\Mo$ or $\Ham\Mo$, 
the Pontriagin 
classes when $\Gg= \Diff(M)$, and the Euler class in the case of $\THha$.
  For example, in the Hamiltonian (but not
the c-Hamiltonian)  case we can use the Chern classes of the
tangent bundle $TM$.  
Each Chern class $c_i(M)$ has a natural extension
$\Tilde c_i$ to $M_\Gg$, namely the $i$th Chern class of the
vertical tangent bundle.  Hence for each multi-index $I: =
(m_1,\dots,m_n)$ there is a class
$$
\mu_{k,I}: = \pi_!\left(\Om^k (\Tilde c_1)^{m_1}\dots (\Tilde
c_n)^{m_n}\right)\in H^*(B\Gg).
$$
As indicated above, initially these classes live on the connected 
group
$\Gg: = \Ham\Mo$.  If $H^1(M) = 0$ they extend to the full 
symplectomorphism
group $\Symp(M)$, but in general only the classes with $k=0$ extend
over this group.  If $M$ is a Riemann surface $\Si$ and $I = m$ then
the classes 
$\ka_m: = \mu_{0,m+1}$ are known as the
 Miller--Morita--Mumford $\ka$--classes.  If $\Si$ has genus $g>1$ 
with 
 orientation class $a$ normalized by $\int_\Si a = 1$,
 then $\Symp(M)$ and 
 $\Hh_a$ are homotopy equivalent to 
 the (orientation preserving) 
 mapping class group $\pi_0(\Symp(\Si))$.  Further
$$
  \Tilde a = \frac 1{2-2g} \bigl(\Tilde c_1 + \frac 
 1{4g-4}\pi^*(\ka_1)\bigr).
  $$
Hence the $\mu_k, k>1,$ together with the single class $\ka_1$ contain 
 the same information as the
$\ka$--classes.

\section{Higher homotopy groups}\label{ss:hh}

We begin with  examples of nontrivial $G=SU(\ell)$ actions for which 
the map $\pi_*(G)\to \pi_*(\Gg)$ is not injective, and then discuss 
an easy way to detect elements in $\pi_*(\Gg)$.   \S\ref{ss:split} 
proves
Proposition~\ref{prop:eval} and \S\ref{ss:gr} concerns flag manifolds.

\subsection{Examples}\label{ss:ex}

\begin{example}[Action of $G=SU(\ell)$ with $r_*= 0$ in dimensions 
    $1$ mod $4$]\label{ex:ex1}\rm
 Let $2\ell\le n+1$ and consider the
action of $SU(\ell)$ on $\cp^n$ given by restricting the standard 
action of $SU(n+1)$
 to the image of the homomorphism
 $$
\rho\co  G = SU(\ell)\to  SU(n+1)\co  A\mapsto (A,\overline A)\in SU(\ell)\times 
 SU(\ell)
 \subset SU(n+1),
 $$
where $\overline A$ denotes the conjugate of $A$.  
Since  conjugation  is an automorphism of $SU(\ell)$ that acts by $-1$ 
on the homotopy
groups in degrees $\equiv 1\pmod 4$, the induced map $\rho_*\co \pi_j(G)\to 
\pi_j(SU(n+1))$ is zero when $j\equiv 1\pmod 4$.  
Hence the image of  $\pi_j(G)$ in $\pi_j(\Ham(\cp^n))$ vanishes for 
these $j$.  Note that the pullback of the  
universal bundle $E_{n+1}\to 
BU(n+1)$ by $B\rho\co  BU(\ell)\to BU(n+1)$ is the sum 
$E_\ell \oplus E_\ell^*\oplus {\B C}^{n+1-2\ell}$.  Thus this example 
reformulates the fact that the odd Chern classes of a bundle of the 
form $E_{\B R}\otimes {\B C}$ vanish (rationally), where $E_{\R}$ denotes the 
real bundle  underlying  $E$.
\end{example}

One can also find representations
$\rho\co  U(\ell)\to U(n+1)$  that kill homotopy in dimensions $4i+3$,
and hence actions of $U(\ell)$ on projective space with the same 
property.

\begin{example}[Action of $SU(4)$ with $r_*= 0$ in dimension $7$]\label{ex:ex2}\rm
If 
$
E\to BU(4)
$
is the universal bundle and $\ga\co  S^8\to 
BU(4))$ generates $\pi_8(BU(4))$, then $c_4(E)(\ga)  = 6$ by 
definition.
 We claim that 
$$
c_4(\La^2(E))(\ga) = -24.
$$
To check this, pull
$c_4(\La^2(E))$  back to $H^*(B\T)$ where $\T$ is the maximal 
$4$--torus in $U(4)$.  If $t_1, \dots, t_4$ denote the obvious 
generators of $H^2(B\T)$ then the pullback of $c_4(\La^2(E))$ to 
$H^8(B\T))$ is
$$
2s_{3,1} + 5 s_{2,2} + 13 s_{2,1,1} + 30 s_{1,1,1,1},
$$
where, for each partition $I = (i_1,\dots,i_r)$ of  the number $4$, $s_I$ denotes 
the symmetric function $\sum t_1^{i_1}\cdots t_r^{i_r}$.
When we express this in terms of the elementary symmetric functions 
$\si_i$ and pull back by $\ga$, the only contribution comes from the 
terms in $\si_4$.  Thus
$$
c_4(\La^2(E))(\ga) = 6(2\cdot 4 + 5\cdot 2 + 13\cdot(-4) + 
30)  = -24.
$$
Therefore the homomorphism $U(4)\to U(22)$ that classifies the sum of
four copies of the rank $4$ bundle $E$ with 
the rank $6$ bundle $\La^2(E)$ kills 
$\pi_7(U(4))$.  A similar statement holds for the corresponding action 
of $SU(4)$ on $\cp^{21}$.
\end{example}

The next result gives a way of detecting the image in $\pi_*(\Gg)$ of 
some of the elements that come from $\pi_*(G_p)$.
Consider the commutative diagram
$$
    \CD
    &&         M\\
    &&         @V jVV \\
    BG_p @>R>> B\Gg_p@> d>> B\Ll\\
    @V q VV    @V\pi VV \\
    BG @>R>>   B\Gg
      \endCD
$$
where $\Gg$ denotes either $\Ham\Mo$ or $\Diff_0(M)$ as appropriate, 
and $\Ll$ is the group formed by the linearized action
of $\Gg_p$  on $T_pM$.  Thus  when $\Gg$ consists of
symplectomorphisms $\Ll = Sp(2n)\simeq U(n)$, while in the smooth
case  $\Ll = GL(2n,\R)\simeq SO(2n)$.  The  homomorphism $d\co  \Gg_p\to \Ll$ is 
given by taking the derivative, so that the composite 
$d\circ j\co  M\to B\Gg_p\to B\Ll$ classifies
the tangent bundle of $M$.  In  degrees
in which this vanishes in 
homotopy we can  use the map $d\circ R$ to detect the image of $R_*$.
For clarity, we state the next lemma in the symplectic case.
Thus $\Gg: = \Ham\Mo$ and $\Ll: = U(n)$ where $2n= \dim M$.

\begin{lemma}\label{le:transpi}  Suppose that the Lie group 
    $G: = SU(\ell)$ acts  on
$M$, and consider the associated map $r\co G\to \Gg$.
Assume that $(d\circ j)_*\co  \pi_k(M)\to \pi_k(B\Ll)$ vanishes, and 
that there is $\ov{\al} \in \pi_k(BG_p)$ such that both $q_*(\ov{\al})\in 
\pi_k(BG)$ and  
$(d\circ R)_*(\ov{\al})\in \pi_k(B\Ll)$ are nonzero.
Then   
 $R_*\co  \pi_k(BG)\oq\to\pi_k(B\Gg)\oq$ is injective.
\end{lemma}
\begin{proof}\, The proof is an easy diagram chase.
    It suffices to show that 
        $$
    (R\circ q)_*(\ov{\al}) = (\pi\circ R)_*(\ov{\al})
    $$
    is nonzero.
    But otherwise $R_*(\ov{\al})$ would be in the kernel of 
    $\pi_*$ and hence in the image of $j_*$.  This would imply
    $d_*(R_*(\ov{\al})) = d_*(j_*(\be))= 0$, which contradicts the 
    assumption that $d_*(R_*(\ov{\al}))\ne 0$.
\end{proof}

We apply this to flag manifolds in Lemma~\ref{le:flags} below.

\subsection{c-split manifolds and the evaluation map}\label{ss:split}

We now prove Proposition~\ref{prop:eval} that detects the image 
of elements 
in $\pi_*(G)$ that map nontrivially under the evaluation map.
 Recall from Lalonde--McDuff~\cite{LM} 
that a fibration $M\to P\to B$ is said to be 
c-split if $\pi_1(B)$ acts trivially on $H^*(M)$
and the Leray--Serre spectral sequence degenerates at the $E_2$
term.  
This is equivalent to the condition
that the inclusion induces an injective map $H_*(M)\to H_*(P)$.
Further a c-symplectic manifold $\Ma$ is said to satisfy the {\em hard
Lefschetz condition} if 
$$
a^k\cup\co  H^{n-k}(M)\to H^{n+k}(M)
$$
is an
isomorphism for all $1\le k\le n$.

\begin{lemma}\label{le:Bl}  Suppose that $\Ma$ satisfies the 
 hard Lefschetz condition and denote by $\Hh$  either $\THha$
or (if  $H^1(M)=0$) the monoid $\Hh_H$ of all homotopy equivalences
that act trivially on rational homology.
Then $M \to M_{\Hh}\to B\Hh$ is
c-split.
\end{lemma}
\begin{proof}  
It is a classical result due to Blanchard~\cite{Bl}
that every fibration over a base $B$ such that 
$\pi_1(B)$ acts trivially on the homology of the fiber
is c-split provided that the class $a$ extends to $P$ and the
fiber $(M,a)$ satisfies the hard
Lefschetz condition.  Hence we need only check that $a$ extends to
$M_{\Hh}$, which is true by 
 Proposition~\ref{prop:ka}.
\end{proof}

{\bf Proof of Proposition~\ref{prop:eval}}\qua  
Consider the fibration $M\to M_{\Hh}\to B\Hh$.  The boundary map $\p\co  
\pi_*(B\Hh)\to \pi_{*-1}(M)$ of its long exact homotopy sequence is 
essentially the same as the evaluation map 
 $\ev_*\co \pi_*(\Hh)\oq\to \pi_*(M)\oq$: more precisely, $\ev =
 \p\circ \tau$
 where $\tau\co \al\to \ov{\al}$ is given by suspension.
Therefore, we 
 must show that if
 a connected  group $G$ acts on a nilpotent manifold $(M,a)$ that satisfies
 the c-splitting condition for $\Hh$, then 
    for every element
$\al\in \pi_{2k-1}(G), k>1,$ with $\ev_*(\al)\ne 
 0$ the element $h(\ov{\al})\in H_*(B\Hh)$ is nonzero.
 To this end, let $ P\to S^{2k}$ be the pullback of
 $M_{\Hh}\to B\Hh$ by  $\ov{\al}$, and 
 consider the  commutative diagram:
 $$
 \CD
M @>=>> M\\
 @VVV  @VVV\\
 P @>>> M_\Hh\\
 @V\pi VV  @V\pi VV\\
 S^{2k}@>>> B\Hh
 \endCD
 $$
 By hypothesis the fibration $M\to M_\Hh\to 
 B\Hh$ is c-split. Hence the elements of
any additive basis $b_0: = 1, b_1,\dots, b_q$ for 
$H^*(M)$ extend to elements $\Tilde b_0: = 1,\dots, \Tilde b_q$ of 
$H^*(M_\Hh)$.  Moreover, if $\deg b_i < 2k$ the restriction of 
$\Tilde b_i$ to $P$, which we denote by $\Tilde b_i^P$,
 is uniquely determined.

Let $e\in H^{2k}(S^{2k})$ be a generator.  Then $\pi^*(e)\ne 0$ since 
the fibration
$P\to S^{2k}$ c-splits.
 We claim that
there is a polynomial $f$ such that the following relation holds in 
$H^*(P)$
for suitable $i_1,\dots,i_s$:
\begin{equation}~\label{eq:e}
\pi^*(e) = f(\Tilde b_{i_1}^P,\dots,\Tilde b_{i_s}^P),\qquad \deg b_{i_j}^P<2k,\;\; 
\forall j.
\end{equation}
To see this, build the  KS model $(\Aa_P,D)$ for the fibration $M\to P\to 
S^{2k}$. 
According to Tralle--Oprea~\cite{TO}, this has the form $(\La(e,t)\otimes \Aa_M, D)$ where
$(\Aa_M,d_M)$ is a minimal model for $H^*(M)$, 
$(\La(e,t),d)$ is a minimal model for $S^{2k}$ (so that
$\deg e = 2k,\deg t = 
4k-1$), and $D$ extends  $d_M$.  
Since $M$ and hence $P$ is nilpotent,
this KS model calculates $H^*(P)$. Since it has
generators in dimensions $\ge 2$ that are 
dual to the rational homotopy, 
it agrees with $\Aa_M$ in dimensions $< 2k-1$. In 
dimension $2k-1$ we may choose generating cochains for $\Aa_M^{2k-1}$ so that 
precisely one of them, say $x$, does not vanish on
the image of $\p\co \pi_{2k}(S^{2k})\to \pi_{2k-1}(M)$.
Then $D(x) = e + r$ where  $r=d_M(x)$ is a product of 
elements in $\Aa_M$ of degrees $< 2k$.  
Since $0= D^2(x) = D(e) + D(r) = D(r)$, $r$ is a cocycle in $\Aa_P$.  Moreover,
because  $\pi^*(e)\ne 0$,
$r$ must represent a nontrivial class in $H^*(P)$.
Therefore if $f$ is the polynomial such that
$f(b_{i_1},\dots, b_{i_s}) = -r$, relation~(\ref{eq:e}) holds.

Now consider the corresponding element
$f(\Tilde b_{i_1},\dots,\Tilde b_{i_s})\in H^*(M_\Hh)$.
Since $\pi\co M_{\Hh}\to B\Hh$ is c-split, we may apply the
Leray--Hirsch theorem.  Therefore
this can be written uniquely in the form
\begin{equation}\label{eq:e1}
f(\Tilde b_{i_1},\dots,\Tilde b_{i_s}) = \sum_{j=0}^q\pi^*(z_j)\Tilde b_j,
\end{equation}
where $z_j\in H^*(B\Hh)$. If we pull this relation back to $P$ we obtain 
the unique expression for $f(\Tilde b_{i_1}^P,\dots,\Tilde b_{i_s}^P)$
in terms of the Leray--Hirsch basis $\Tilde b_j^P$ for $H^*(P)$.
Comparing with (\ref{eq:e}), we find that $z_0\in H^{2k}(B\Hh)$ must 
extend $e$.  Thus 
$\ov{\al}$ has nonzero image in 
$H_*(B\Hh)$.
\qed

\begin{rmk}\label{rmk:xx}\rm (i)\qua  Equation~(\ref{eq:e}) implies that
    $\pi_!(f(\Tilde b_{i_1}^P,\dots,\Tilde b_{i_s}^P)) = e$.   In 
    fact the fiber integral simply picks out the coefficient of 
    $b_0: = 1$ 
    in the 
    Leray--Hirsch decomposition for an element in $H^*(M_{\Hh})$. 
    \ms
    
    (ii)\qua 
    The c-splitting hypothesis in
    Proposition~\ref{prop:eval} is satisfied when $M$
is a c-symplectic manifold with a transitive action of
a connected and simply connected compact 
Lie group $G$.  For if $\om$ is
a $G$--invariant representative 
for the c-symplectic class  $a\in H^2(M)$, then
$\om^n$ has constant rank  and hence does not vanish.
Thus $\om$ is symplectic.   Our assumptions  imply that $G$ is 
semisimple and that $\pi_1(M)=0$.  Therefore  $M$ has a K\"ahler 
structure with K\"ahler form $\om$ 
by a theorem of Borel that is stated for example in 
Tralle--Oprea~\cite[Ch~5,Thm~2.1]{TO}.
Therefore the claim follows from Lemma~\ref{le:Bl}.
\end{rmk}

We now give a version of 
 Proposition~\ref{prop:eval} in which the c-splitting condition is 
 relaxed.   All we need is that the cohomology classes in $M$ occuring 
 in the relation (\ref{eq:e}) extend to $H^*(M_{\Hh})$.  In 
 particular, if in the statement below $H^*(X)$ is generated by the 
 Chern classes of $M$ then we can forget about the class $[\om]$ and 
 take $\Gg: = \Symp(M)$
 even if $H^1(M)\ne 0$.

\begin{prop}\label{prop:eval2}  Suppose that a simply connected group
    $G$ acts on $(M,\om)$ in 
such a way that the subring of $H^*(M)$ generated by $[\om]$ and the Chern 
classes of $TM$ surjects onto $H^*(X)$, where $X: = G/G_p$ is some 
orbit of $G$.  
Let $\Gg$ be a subgroup of $\Symp(M)$ that contains $G$ and is such that 
$[\om]$ extends to $H^*(M_\Gg)$.  Then the map
$$
h\circ R_*\co  \pi_*(BG)\to H_*(B\Gg)
$$
is nonzero on all elements $\ov{\al}\in \pi_{2k}(BG)$ such that
$\ev_*(\al) \in \pi_{2k-1}(M)$ is nonzero.
\end{prop}
\begin{proof}  Suppose given $\al \in \pi_{2k-1}G$
    such that $\ev_*(\al)\ne 0$ in $\pi_{2k-1}(M)$ and hence
   also in $\pi_{2k-1}(X)$.  Consider the corresponding diagram:
$$
 \CD
X @>>> M @>=>> M\\
 @VVV  @VVV  @VVV\\
Q @>>> P @>>> M_\Gg\\
 @V\pi VV @V\pi VV  @V\pi VV\\
 S^{2k}@>=>> S^{2k}@>>> B\Gg
 \endCD
$$    
   Our hypotheses ensure that 
 there is a basis $b_0: = 1, b_1,\dots, b_N$ for $H^*(X)$ consisting 
 of elements that 
 extend to $H^*(M_\Gg)$.
 Hence the fibration  $\pi\co Q\to S^{2k}$ is c-split
 and the result follows by arguing
 as in the proof of
 Proposition~\ref{prop:eval}.
\end{proof}

\subsection{Flag manifolds}\label{ss:gr}

Consider the general flag manifold
$$
M(m_1,\dots,m_k): = U(\ell)/U(m_1)\times\cdots\times U(m_k),\qquad
m_1\ge\dots\ge m_k,
$$
where $\ell = \sum m_i$.  We shall denote 
$$
G': = U(\ell),\quad  G_p': = U(m_1)\times\cdots\times U(m_k),\quad
G: = SU(\ell), \quad G_p: = G_p'\cap G.
$$
If $E_i\to U(m_i)$ denotes the universal bundle, then 
 the vertical tangent bundle of the fibration
$$
\CD
M @>j>> BG_p@>q>> BG
\endCD
$$
is the pullback 
from $BG_p'$ to $BG_p$ of the bundle
$$
\bigoplus_{1\le i<j\le k}\Hom(E_i,E_j) = \bigoplus_{1\le i<j\le 
k} E_i^*\otimes E_j
$$
Further,  the cohomology ring
of $M$ is generated by the Chern classes of 
the $E_i$ with defining relations coming from the fact that
the  restriction  of $\oplus_i E_i$ to $M$ is  trivial.  In 
particular,  in the case of a Grassmannian (ie $k=2$) the cohomology 
of $M$ is generated by the 
Chern classes of its tangent bundle
so that $\Symp\Mo$ acts trivially on $H^*(M)$.

The first part of the next lemma illustrates the use of 
Lemma~\ref{le:transpi}, though  in this case stronger conclusions
may be obtained by other methods (see Proposition~\ref{prop:flags2}.)
Recall that $\Hh_H$ 
denotes  the space of homotopy equivalences that act 
 trivially on $H^*(M)$.

\begin{lemma}\label{le:flags}
Let $M$ be the
flag manifold $M(m_1,\dots,m_k)$ and set $G: = SU(\ell)$ where
$\ell = \sum m_i$ and $m_1\ge\dots\ge m_k$. 
Then:

{\rm (i)}\qua $r_*\co \pi_{2i-1}(G)\to \pi_{2i-1}(\Symp\Mo)$ is injective for
$m_2<i\le m_1$. 
\smallskip

{\rm (ii)}\qua    $h\circ R_*\co \pi_{2i}(BG)\to 
H_{2i}(B\Hh_H)$ is injective for  $i = 2$ and $m_1<i\le \ell$.
\end{lemma}
\begin{proof}\,   The existence of the fibration
    $$
M\to BG_p\to BG
    $$
implies that the generators of  $\pi_*(M)$ divide into two groups.
Those in dimensions $2i, 1\le i\le m_2,$ map to the elements in the 
kernel of  $\pi_*(BG_p)\to\pi_*(BG)$, while the 
generators in the image of $\ev_*$ come from $\pi_{*+1}(BG)$ and
lie in odd dimensions $2i-1$ for 
$m_1< i\le \ell$. 
Therefore because 
$\pi_{2i}(M) = 0$ in the 
range $m_2<i\le m_1$, the
 first statement will follow from Lemma~\ref{le:transpi} if we 
show that $(d\circ r)_*\co  \pi_{2i}(BG_p) \to \pi_{2i}(B\Ll)$ is 
surjective for $m_2< i \le m_1$.  
To prove this it suffices to check that the Chern classes of the 
vertical tangent bundle of the fibration $M\to BG_p\to BG$
are nonzero on  $\pi_{2i}(BG_p)$  for $i$ in this range.
This holds because the restriction of this bundle to $BSU(m_1)$ is 
simply a sum of copies of 
$E_1^*$. This proves (i).

Since $\ev_*\co \pi_*(G)\to \pi_*(M)$ is injective in dimensions $*>2m_1-1$,
(ii) follows immediately from
Propositions~\ref{prop:nonzero}
and ~\ref{prop:eval} and Lemma~\ref{le:Bl}.
\end{proof}

We next show that the map 
$h\circ R_*\co \pi_*(BG)\to 
 H_*(B\Hh_H)$ is always injective for flag manifolds.  Our argument 
imitates the proof for $\cp^n$,
 considering the coefficients of the higher degree 
elements $\Tilde b_j, j> 0,$ in (\ref{eq:e1}).

\begin{prop}\label{prop:flags2}  Let $M$ be the flag manifold
    $M(m_1,\dots,m_k)$  and set $G: = SU(\ell)$ where
$\ell = \sum m_i$ and $m_1\ge\dots\ge m_k$. 
 Then the inclusion 
$
BSU(\ell)\to 
B\Hh_H
$
induces an injection on $H_*$.
\end{prop}
\begin{proof} Let us first consider the case of a Grassmannian
    $M(m,k)$ where $\ell = m+k, m\ge k,$ and set $G: = SU(\ell)$.  
    It is convenient to work with 
    cohomology.  Therefore we aim to show that the Chern classes 
    $c_i, i=2,\dots,\ell$ in $H^*(BG)$ extend to classes in $H^*(B\Hh_H)$.
    
As remarked above, 
the ring $ H^*\bigl(B(U(m)\times U(k))\bigr)$ is freely generated by elements
$x_1,\dots,x_m$ (the pullbacks of the Chern classes 
 of the universal bundle $E: = E_1\to BU(m)$) and 
$y_1,\dots,y_k$ (the  pullbacks of the 
Chern classes  
of the universal 
bundle $F : = E_2\to BU(k)$).  Thus 
$H^*(BG_p) = H^*\bigl(BS(U(m)\times U(k))\bigr)$ is the quotient of this free 
ring by the relation $x_1+y_1=0$.  

Consider the fibration $\pi\co  
BG_p\to BG$.
Since $E\oplus F$ is the pullback of the universal bundle over
$BG: = BSU(\ell)$  its Chern classes are the pullbacks
$Q_i: = \pi^*(c_i)$ of the Chern classes in $H^*(BG)$. 
Thus $Q_0: = 1$ and $Q_1 = 0$.
Taking the total Chern class of $E\oplus F$ we find:
\begin{equation}\label{Q}
(1+x_1+\dots+x_m)(1+y_1+\dots+y_k)=1+Q_1+\dots+Q_{m+k}.
\end{equation}
Since the restriction of $E\oplus F$ to $M$ is trivial, the above 
identity gives $m+k$ relations among the restrictions of the 
$x_i, y_j$ to $H^*(M)$.  The first $m$ of these 
should be interpreted as defining  the $x_i$'s 
in terms of the generators $y_i$ of $H^*(M)$, while the 
equations  $Q_{m+i} = 0$ give the
relations  in $H^*(M)$.  In particular, there 
are no relations among the monomials in the 
$y_i$ in degrees $\le 2m$.
Hence we may choose  an additive basis $b_0: = 1, b_1,\dots, b_N$ 
for $H^*(M)$ whose elements in degrees $\le 2m$ consist of all 
the monomials in the 
$y_i$.  
In these degrees we may therefore write $b_\nu = y_{I_\nu}$,
and will extend $b_\nu\in H^*(M)$ to $H^*(M_G)$ by identifying it with 
$y_{I_\nu}\in H^*(M_G)$.
As in Proposition~\ref{prop:eval}, the $b_\nu$ also extend to
elements $\Tilde b_\nu\in H^*(M_{\Hh})$.  We denote by $\Tilde b_\nu^G$ the 
restriction of $\Tilde b_\nu$ to $M_G = BG_p$. These two 
extensions $\Tilde b_\nu^G$ and $b_\nu$ of $b_{\nu}\in H^*(M)$ to 
$H^*(M_G)$ need not agree, but they do agree modulo the 
ideal $\bigl\langle Q_2,\dots,Q_m\bigr\rangle$ generated by the 
elements of $\pi^*(BG)$.  In other words, if we identify
 $H^*(BG_p)$ with the free algebra generated by the $y_j$ 
and $Q_i, 2\le i\le m$, we have
$$
b_\nu \in  \Tilde b_\nu^G + \bigl\langle Q_2,\dots,Q_m\bigr\rangle.
$$

 Formally inverting  $1+y_1+\dots +y_k$ in equation~(\ref{Q}) we obtain
$$
1+x_1+\dots +x_m = (1+Q_1+\dots+Q_{m+k})(1+f_1(y_1,\dots,y_k)+\dots),
$$
where the $f_i$ are the homogeneous terms of degree $2i$ in
$(1+y_1+\dots+y_k)^{-1}$.
The terms of degree $2m+2$ give the relation
$$
-Q_{m+1}= \sum_{i=1}^{m}Q_{m-i+1}f_i(y_1,\dots,y_k).
$$
Notice also that the coefficient of $y_1^i$ in $f_i$ equals $(-1)^i$.
Since the polynomials $f_i$ have degree $2i\le 2m$  they are
sums $\sum \al_{\nu} b_{\nu}$ of the basis monomials $b_\nu = y_{I_\nu}$.
Therefore we have
\begin{eqnarray*}
\pi^*(c_{m+1})=
Q_{m+1} &=&- \sum_{i=1}^{m}Q_{m-i+1}f_i(y_1,\dots,y_k)\\
& = & -\sum_{i=1}^{m}Q_{m-i+1}\left(\sum_\nu\al_{\nu}\,b_{\nu}\right)\\
& = & \sum_\nu \pi^*(e_\nu) \,\,\widetilde b^G_\nu
\end{eqnarray*}
where each $e_\nu\in H^*(BG)$ is a 
polynomial in the Chern classes with a nonzero linear term.
Further, if $ b_{\nu_i}$  denotes 
the basis element $y_1^{m-i}$ of $H^*(M)$ for $1< i\le m,$  then
the coefficient  of $\Tilde b_{\nu_i}^G$ in the above expression
has the form $\mu c_i +$ decomposables, for some nonzero number $\mu$.

It follows from Proposition~\ref{prop:eval} that there is a class
$u\in H^{2m+2}(B\Hh)$ whose restriction to $BG$ does not vanish on 
$\pi_{2m+2}(BG)$.  Hence we may choose $u$ so that it restricts to
$c_{m+1} + p(c)\in H^{2m+2}(BG)$ where $p(c)$ is some polynomial in 
the $c_i, 2\le i\le m$.     
Thus $\pi^*(u)$ restricts to $Q_{m+1} + \pi^*(p(c))$ 
where $\pi^*(p(c))$ is a polynomial in the $Q_i, i\le m$.  Since
 $\pi^*(u)$ vanishes on the fiber $M$ and there is a unique relation 
 in $H^{2m+2}(M)$ namely $Q_{m+1}=0$, $\pi^*(p(c))$ must be a 
 multiple of $Q_{m+1}$.  But $H^*(BG_p) = H^*(M_{G})$ is freely 
 generated by the $y_i$ and $Q_j, j>1$, and so this is possible only if 
$\pi^*(p(c)) = 0$.
Hence $\pi^*(u)$ extends $Q_{m+1}$ and the
argument may be completed as before.
The Leray--Hirsch theorem implies that
$\pi^*(u)$ may be written uniquely as
$$
\pi^*(u) = \sum_\nu \pi^*(u_\nu) 
\Tilde b_\nu \;\;\in H^*(M_\Hh),
$$
where $u_\nu\in H^*(B\Hh)$.  Comparing with the expression previously 
found for $Q_{m+1}$ we see that a multiple of $u_{\nu_i}$ extends 
a class which equals $c_i$ modulo products of $c_j, j<i$.  An easy 
inductive argument now shows that each $c_i$ must extend to $B\Hh$.
This completes the proof for Grassmannians.

The proof for the flag manifold $M(m_1,\dots,m_k),k>2,$ is very similar.
We denote by $x_i$ the Chern classes of the 
universal bundle $E_1\to BU(m_1)$ and by  $y_{\al i}$ the Chern classes of
$E_\al\to BU(m_\al)$ for $2\le\al\le m$. Further let $y_i$ be the
Chern classes of the sum $F:=E_2\oplus\dots\oplus E_m$.
Since $E_1\oplus F$ is pulled back from $BU(\ell)$ it is trivial on 
$M$.  Further the equation~(\ref{Q}) 
 holds as before.  It should be interpreted as first
defining the $x_i$ in terms of the generators $y_{\al j}$ for $H^*(M)$
and then giving the relations in $H^*(M)$.
The rest of the argument goes through without essential change.
\end{proof}

\section{The evaluation map and Whitehead products}\label{sec:eval}

In this section we consider a nontrivial
 $S^1$ action $\la$ that is inessential in $\Gg$.
 For simplicity, we suppose throughout this section 
 that $M$ is simply connected, though versions of the first lemmas  
 extend to the general case.
 We shall suppose either that we are in the symplectic category so
 that $\Gg = \Ham\Mo$ or that $\Ma$ is c-symplectic, the action is smooth,
 and $\Gg: = \Diff_0(M)$.
Let $\rho\in \pi_3(\Gg)$ be  the nonzero element
constructed in Theorem~\ref{thm:s0}.  Our first aim is to understand
what it means for $\ev_*(\rho)\in \pi_3(M)$ to vanish.

We denote by $\p$ the boundary map
in the long exact sequence of the fibration
$$
\CD
\Gg_p@>>> \Gg @>\ev>>M.
\endCD
$$

\begin{lemma} \label{le:al} Let $M,\Gg$ be as above
    and suppose that $\la$ is inessential in $\Gg$
    but essential in $\Gg_p$.  Then
    there is a unique
    $\al\in \pi_2(M)\oq$ such that $\p\al = \la$.  Moreover,
    $h(\al)\in H_2(M)$ is nonzero.
\end{lemma}
\begin{proof}  The first statement holds because
the evaluation map
$$
\ev_*\co \pi_2(\Gg)\oq\to
\pi_2(M)\oq
$$
is zero so that  $\p$ is injective. 
The second is an immediate 
consequence of the Hurewicz theorem: 
because $\pi_1(M)=0$ the map
$\pi_2(M)\oq\to H_2(M)$ is an isomorphism.
\end{proof}

Note that the element
 $\al$ may depend on the choice of $p$.   For example the
loop in $\Ham(\cp^2,\om)$ given by $[z_0:z_1:z_2]\mapsto
[e^{2\pi it}z_0:e^{-2\pi it}z_1:z_2]$ is nullhomotopic in $\Gg_p$
when $p = [0:0:1]$ but is essential when $p = [1:0:0].$  Therefore if 
we work with the corresponding diagonal circle action on 
$M=\cp^2\times\cp^2$ we can find points $p_1, p_2$ such that 
$\la$ is essential in $\Gg_{p_i}$ for $i = 1,2$ but the elements 
$\al_1,\al_2$ are different.

\begin{lemma} \label{le:alal} Let $\la$, $\al$ and $\rho: =
    \{\la,\la\}$ be as in Lemma~\ref{le:al}. Then  $\ev_*(\rho) $ is the 
    Whitehead product $ [\al,\al]$.
    \end{lemma}
\begin{proof}  Consider the fibration $\Om M \to \Tilde{\Gg}_p\to 
\Gg$,
    where $\Om M$ is the based loop space of $M$ and  $\Tilde{\Gg}_p$
    is the space of all pairs $(h,\ga)$ where $h\in \Gg$ and
    $\ga$ is a path in $M$
    from the base point $p$ to $h(p)$.  The choice of a contraction
    $\Tilde\la$ of $\la$ in
    $\Gg$ determines a homotopy from $\la$ to a loop
    $\ell\co s\mapsto \ell(s)$ in $\Om M$: if $\Tilde\la$ is given by a 
map
\begin{eqnarray*}
 &&   [0,1]\times [0,1]\to \Gg, \qquad (s,\nu) \mapsto 
    \Tilde\la(s,\nu),\\ 
 &&   \Tilde\la(s,0) = \Tilde\la(0,\nu) = \Tilde\la(1,\nu) = id,\quad
    \Tilde\la(s,1)=\la(s),
\end{eqnarray*}
  then for each $s\in [0,1]$, $\ell(s)$ is the loop $\nu\mapsto
  \Tilde\la(s,\nu)(p)$.

  It suffices to show that the image of $\rho$ under the boundary map
  $\p\co \pi_3(\Gg)\to \pi_2(\Om M)$ is the Samelson product
  $\left<\ell,\ell\right>$ since this is the desuspension of
  $[\al,\al]$.  But it follows from the definition of
  $\rho$ as a map $D^2\times S^1/\!\sim\;\;\Rar \;\Gg$ that $\p\rho$ 
is
  represented by the map
  $$
  S^1\times S^1/\!\sim\;\;\Rar\;\; \Om M,\quad (s,t)\mapsto
  \left<\Tilde\la(s,\cdot),\la(t)\right> (p),
  $$
  where $(s,t)\in S^1\equiv \R/\Z$.  But this is homotopic to the map
  $
  \left<\ell,\ell\right>\co  (s,t)\mapsto\ga(s,t)$ where 
  $$
 \ga(s,t)(\nu): =  \left<\Tilde\la(s,\nu),\Tilde\la(t,\nu)\right>(p)
  $$
  via the homotopy $\ga_r(s,t)(\nu): =
  \left<\Tilde\la(s,\nu),\Tilde\la(t,r + \nu(1-r))\right>(p)$.
  \end{proof}

  A similar argument shows that if the product
  $\{f,f'\}$ is as defined in 
Proposition~\ref{prop:uniq} then $\ev_*(\{f,f'\})  = [\al,\al']$, 
provided that 
  $f,f'$ both vanish in $\Gg$ so that they 
  correspond to elements $\al,\al'\in \pi_*(M)$.  
  If only $f$ vanishes
  in $\Gg$, $\ev_*(\{f,f'\})$ need not be a Whitehead product.
 For example, in the case of the action of
  $SU(n+1)$ on $\cp^n$ we have (in the notation of Lemma~\ref{le:un})
  $\ev_*(\{\la_1,\la_3\}) = 0$, but
  $\ev_*(\{\la_1,\la_{2n-1}\})\ne 0$.   In this particular example,
  $\ev_*(\{\la_1,\la_{2n-1}\})$ is a higher order Whitehead product.
  We now investigate the extent to which this generalizes.
  To get clean statements we need to assume that the homotopy of
  $M$ has some of the characteristics of $\cp^n$.  For the first 
  result we only need $M$ to be simply connected, but later on need
  more assumptions.

\begin{prop}\label{prop:al2} 
Let  $\la$ be a nontrivial
    $S^1$--action on a simply connected c-symplectic manifold $(M,a)$.
    Suppose that there
    is $\al\in \pi_2(M)$ 
    such that $\p\al = \la\in \pi_1(\Gg_p)$.  Then
    $\ev_*(\rho) \ne 0$ if and only if  any element $c\in H^2(M)$ such 
that
    $c(\al)\ne 0$ has the property that $c^2 = 0$ modulo the ideal
    $\Ii$  in $H^*(M)$ generated by
    the kernel of $\al^*\co  H^2(M)\to H^2(S^2)$.
\end{prop}
\begin{proof}\,
If $\ev_*(\rho) = 0$ then Lemma~\ref{le:alal} implies
that the map $\al\vee\al\co S^2\vee S^2\to M$ extends to $S^2\times S^2$.
Hence if $c\in H^2(M)$ is such that $c(\al)\ne 0$ its square
$c^2$ has nonzero pullback
to $S^2\times S^2$ and hence is nonzero modulo $\Ii$.

The converse follows from from Lemma~\ref{le:alal} and minimal
    model theory.  When building a minimal model $\Aa: = 
(A^k,d)_{k\ge 0}$
    for $H^*(M)$ we may
    choose a basis  $c_0,\dots c_k$ for the $2$--cochains $A^2$
    such that
    $c_i(\al)=0, i>0$.  The differential
    $d\co  A^3\to (A^2)^2\subset A^4$ is dual to
    the Whitehead product $\pi_2(M)\times \pi_2(M)\to \pi_3(M)$.
    We may choose a basis for $\pi_3(M)$ whose first element is
    $\ev_*(\rho)$.  Then if $r
    \in A^3$ is the first element in the dual basis,
     $dr = \sum\mu_{ij}c_i c_j$ where $\mu_{00} = 1$.  Hence $dr =
     (c_0)^2$ modulo $\Ii$.
  \end{proof}

We now explore what happens when $ev_*(\rho)=0$. This means that
the Whitehead product $[\al,\al]$ is zero. Thus we can consider the
higher order Whitehead product $[\al,\al,\al]\subset \pi_5(M)$
which is defined as follows (see \cite{AA} for details). 
Let $W\subset S^2\times S^2\times S^2$ 
denote the fat wedge, that is it consists of triples with at least
one coordinate at the base point. 
Then $[\al,\al,\al]:= \{f_*(u)\}\subset \pi_5(M)$, where
$f\co W\to M$ ranges over the set of all extensions of $\al\vee\al\vee\al$ and
$u\in \pi_5(W)$ is a generator. Since $[\al,\al]=0$ 
this set is nonempty.  Moreover, because $W$ consists of the wedge of 
three copies of 
$(S^2\times S^2)\vee S^2$, $f$ is determined by the 
way in which $\al\vee\al$ is extended to $S^2\times S^2$, which can 
vary by an element $\be\in \pi_4(M)$.  Therefore 
$[\al,\al,\al]$ is a coset of the subgroup
 $$
 H: = \{[\al,\be]: \be \in
\pi_4(M)\}\subset \pi_5(M).
$$
We say that $[\al,\al,\al]$ is nonzero if this coset
does not contain the zero element.  Hence $[\al,\al,\al]=0$ if and 
only if some extension $f\co W\to M$ of $\al\vee\al\vee\al$
extends further to $S^2\times S^2\times S^2$.

In order to interpret the vanishing of $\ev_*(\rho)=0$ in cohomological 
terms (rather than in terms of the minimal model)
we need to make some further simplifying assumptions about the homotopy 
type of $M$.  We assume below that $\pi_3(M)= 0$, so that 
$[\al,\al] = 0$ {\it a fortiori}.  Note that this hypothesis is
satisfied by all generalized
flag manifolds $M(m_1,\dots,m_k)$ with $m_1> 1$.

\begin{proposition}\label{prop:wh}
Let $\pi_1(M)=\pi_3(M)=0$ and
$\la\co S^1\to \Gg_p$ be a nontrivial action as in Proposition
\ref{prop:al2}. Suppose that there is $\al\in \pi_2(M)$ such that
$\partial \al = \la \in \pi_1(\Gg_p)$. Consider $c\in H^2(M)$ such
that $c(\al)\neq 0$.
Then the higher Whitehead product $[\al,\al,\al]$ is nonzero
if and only if $c^3\in \Ii$, where $\Ii$ is as in 
Proposition~\ref{prop:al2}.
\end{proposition}

\begin{proof}
One direction is trivial. Namely,
if $[\al,\al,\al]=0$ then $\al\vee\al\vee\al$ extends to
 $S^2\times S^2\times S^2$ and so $c^3 \notin \Ii$.

Now assume that $[\al,\al,\al]\neq 0$. Chose some extension
of $\al\vee\al$ to $S^2\times S^2$ and using this define
$f\co W\to M$ in a symmetric way.
Denote $\ga_0: = f_*(u)$. By hypothesis
 the span in $\pi_5(M)$ of the
Whitehead products $[\al,\pi_4(M)]$ does not include $\ga_0$.

The minimal models of the spaces under consideration have the
following forms:

\begin{itemize}
\item
$\mathcal A(W) = \La (x_1,x_2,x_3;y_1,y_2,y_3;w;...),$

where $\deg x_i =2$, $\deg y_i = 3$, $\deg w =5$ and
the nontrivial differentials are $dy_i=x_i^2$,
$dw = x_1x_2x_3$ and 
enough others in higher degrees to cancel out the cohomology.
Note that this minimal model is infinitely generated.

\item
$\mathcal A(M) = \La (c_0,c_1,...,c_k;v_1,...,v_\ell ;z...)$,

where $\deg c_i=2$, $\deg v_i=4$, $\deg z = 5$ and
the first nontrivial differential is given by
$dz = \sum a_{ijk}c_ic_jc_k + \sum a_{\ell m}c_\ell v_m$,
for some $a_{ijk},a_{\ell m}\in \B Q$. Moreover,
$c_0(\al)=1$ and $c_i(\al)=0$ for $i>0$.
\end{itemize}

According to the above discussion about
$[\al,\pi_4(M)]$  we can choose $z\in A^5(M)$  so that
$z$ vanishes on all the elements $[\al,\be]$ for $\be \in \pi_4(M)$
but $z(\ga_0)= 1$. 
>From this we draw two conclusions.
Firstly,  because the quadratic part 
of $dz$ is dual to the
Whitehead product,
$dz$ includes no terms of the form $c_0 v_m$. In other words $a_{0m}= 0$ 
for all $m$ and $dz = \sum a_{ijk}c_ic_jc_k$.  
Secondly, because
$f^*(z) = w + \sum b_{ij}x_i y_j$ and $f^*(c_i) = 0, i>0$,
\begin{eqnarray*}
0&\ne &x_1x_2x_3 + \sum b_{ij} x_ix_j^2\\
& = & d(f^*(z))\;\;=\;\; f^*(dz) 
= a_{000}f^*(c_0^3).
\end{eqnarray*}
Therefore $a_{000}\neq 0$, and
the relation in cohomology given by setting $dz$ to zero
says that
$[c_0]^3\in \Ii$.
\end{proof}

Our final result concerns the question of whether
$[\al,\al,\al]\cap \im \ev_*\neq \emptyset$, that is
whether the Whitehead product contains elements from the
image of the evaluation map.  Again, we need to strengthen the 
hypothesis that $[\al,\al]=0$.  The latter implies that 
 $\rho$  lifts to an element
 $\rho_1\in \pi_3(\Gg_p)$, 
 and we now assume that this lift can be chosen 
so that the Samelson product $\langle\la,\rho_1\rangle$ vanishes  in
$\Gg_p$.
By Lemma~\ref{le:un}, this will hold if, for example, 
$\la$ is an essential  circle in  $U(2)\subset G_p$
and the map $U(2)\to G$ factors through $SU(2)$.
Again, this hypothesis holds for flag manifolds $M(m_1,\dots,m_k)$ 
with $m_1>1$, and in particular for $\cp^n$.

\begin{prop}\label{prop:wh2}  Let $\la$ be a nontrivial
    $S^1$--action on a simply connected c-symplectic manifold
    $\Ma$, and let $\Gg: = \Ham\Ma$.  Suppose  that there
    is $\al\in \pi_2(M)$  such
     that $\p\al = \la\in \pi_1(\Gg_p)$ and $[\al,\al] = 0$.
    Suppose further that
    $\rho$ has a lift $\rho_1$ to 
    $\Gg_p$ such that $\langle\la,\rho_1\rangle = 0$ in $\Gg_p$.
Then if $[\al,\al,\al]\ne 0$ 
in $\pi_5(M)\oq$ there is a nonzero element in the 
intersection of $[\al,\al,\al]$ with the
    image of $\ev_*\co \pi_5(\Gg)\to \pi_5(M)$.
\end{prop}

\begin{proof}
Consider the commutative diagram:
$$
\CD
S^2\vee S^2 @>\al\vee\al >> M\\
@VVV   @V j VV\\
S^2\times S^2 @>\La\times \La>> B\Gg_p
\endCD
$$
By assumption, the map $\al\vee\al$ has some extension $\phi$
to $S^2\times S^2$.  However,  because $\pi\circ (\La\times \La)$ is 
nontrivial (where $\pi\co B\Gg_p\to B\Gg$) we cannot choose $\phi$ to 
make this diagram commute.  On the other hand, because 
$\pi\circ (\La\times \La)$ is the composite of the quotient map 
$S^2\times S^2\to S^4$ with
 a representative of  $\ov{\rho}$,
we can adjust $\La\times \La$ on the top cell of $S^2\times S^2$
by a lift $\tau$ of $-\ov{\rho}$
to a map $f: =\La\times \La \#\tau\co  S^2\times S^2\to B\Gg_p$ that does 
lift to $M$.  Thus we can arrange that $f = j\circ\phi$
for suitable $\phi$, ie, that the following diagram commutes:
$$
\CD
S^2\times S^2 @>\phi >> M\\
@VVV   @V j VV\\
S^2\times S^2 @>f>> B\Gg_p
\endCD
$$
Our assumptions imply that we may choose $\tau$ so that 
in addition $[\La,\tau] = 0$ in $B\Gg_p$. 
Consider the commutative diagram:
$$
\CD
(S^2\times S^2\times S^2)_4 @>\Phi >> M\\
@VVV   @V j VV\\
(S^2\times S^2\times S^2)_4 @>F >> B\Gg_p
\endCD
$$
where $\Phi$ is given by $\phi$ and $F$ is given by $f$. Since 
$\La\vee\La\vee\La$ extends to the product map $\La\times \La\times 
\La$, the
obstruction to extending $F$ to the product
$S^2\times S^2\times S^2$ is $3[\La,\tau]$ and so vanishes.  
Therefore the obstruction 
to extending $\Phi$ is an element in $[\al,\al,\al]$ that has zero 
image in $\pi_*(B\Gg_p)$ and hence lies 
in the image of $\ev_*$.
\end{proof}

If $[\al,\al,\al] = 0$ one could explore 
yet higher order products.
However, to make sense of the results one would need further 
assumptions on $M$ and the behavior of $\la$.



\begin{thebibliography}

\bibitem{AM}  {\bf M Abreu}, {\bf D McDuff}, {\it Topology of symplectomorphism
      groups of rational ruled surfaces},
      J. Amer. Math. Soc. 13 (2000) 971--1009
  \MR{1775741}

\bibitem{All} {\bf C Allday}, {\it Examples of circle actions on
      symplectic spaces}, from: ``Homotopy and Geometry (Warsaw
      1997)'', Banach Center Publications 45, Institute of
      Mathematics, Polish Acad. Sci. Warszawa (1998)
      \MR{1679851}

\bibitem{AA} {\bf P Andrews}, {\bf M Arkowitz}, 
       {\it Sullivan's minimal models and higher order Whitehead products}, 
      Canad. J. Math. 30 (1978) 961--982
  \MR{0506254}

\bibitem{An} {\bf S Anjos}, {\it Homotopy type of symplectomorphism groups of 
$S\sp 2\times S\sp 2$}, \gtref6{2002}7{195}{218} 
  \MR{1914568}

\bibitem{AG} {\bf S Anjos}, {\bf G Granja}, {\it Homotopy
decomposition of a group of symplectomorphisms of $S\sp 2\times S\sp
2$}, Topology 43 (2004) 599--618 \MR{2041632}


\bibitem{Bl}  {\bf A Blanchard},  {\it Sur les vari\'et\'es analytiques
      complexes}, Ann. Sci. Ec. Norm. Sup. 73 (1956)
      157--202
  \MR{0087184}

      
\bibitem{GalK}  {\bf S Gal}, {\bf J K\c edra}, {\it Symplectic configurations},
                in preparation

\bibitem{GGK}  {\bf V Ginzburg}, {\bf V Guillemin}, {\bf Y Karshon},  
      {\it Moment maps, cobordisms and Hamiltonian group actions},
      Mathematical Surveys and Monographs 98, AMS (2002)
  \MR{1929136}


\bibitem{GOTT} {\bf D Gottlieb}, 
        {\it Evaluation subgroups of homotopy groups},
        Amer. J. Math. 91 (1969) 729--756 
  \MR{0275424}

\bibitem{G}
{\bf M Gromov},  {\it Pseudoholomorphic curves in symplectic manifolds}, 
Invent. Math.  82  (1985) 307--347 
  \MR{0809718}

\bibitem{GLS} {\bf V Guillemin}, {\bf E Lerman}, {\bf S Sternberg}, 
       {\it Symplectic fibrations and multiplicity diagrams\/}, Cambridge 
       University Press (1999)
  \MR{1414677}

\bibitem{gs}  {\bf V Guillemin}, {\bf S Sternberg},
       {\it Supersymmetry and Equivariant de Rham Theory},
       Springer (1999)
  \MR{1689252}

\bibitem{JK} {\bf T Januszkiewicz}, {\bf J K\c edra}, {\it
Characteristic classes of smooth fibrations}, \arxiv{math.SG/0209288}

       
\bibitem{LM} {\bf F Lalonde}, {\bf D McDuff}, {\it Symplectic structures on
      fiber bundles}, Topology 42 (2003)
      309--347
  \MR{1941438}

\bibitem{Mcf}   {\bf D McDuff},   {\it Symplectic diffeomorphisms and the flux
          homomorphism},  Invent. Math. 77 (1984) 353--66
  \MR{0752824}

\bibitem{Mcox}  {\bf D McDuff}, {\it A survey of topological properties of 
       groups of symplectomorphisms},  from:     
       ``Topology, Geometry and Quantum Field Theory, 
       Proceedings of 2002 Symposium
       in honor of G\,B Segal'', (U\,L Tillmann, editor), LMS
       Lecture Notes series 308, Cambridge Univ. 
       Press (2004) 173--193
  \MR{2079375}

\bibitem{Me}  {\bf D McDuff}, {\it Enlarging the Hamiltonian group}, preprint 
       (2004)

\bibitem{MS} {\bf D McDuff}, {\bf D Salamon}, {\it Introduction to
       Symplectic Topology}, 2nd edition, Oxford Univ. Press, Oxford,
       (1998) \MR{1698616}


\bibitem{MT1} {\bf D McDuff}, {\bf S Tolman},
      {\it Topological properties of Hamiltonian circle actions}, 
      \arxiv{math.SG/0404338}
      
\bibitem{MT2} {\bf D McDuff}, {\bf S Tolman}, in preparation
      
\bibitem{Ono} {\bf K Ono}, {\it Floer--Novikov cohomology and the flux 
        conjecture}, preprint (2004)

\bibitem{Rez} {\bf A\,G Reznikov}, {\it Characteristic classes in symplectic
      topology},  Selecta Math 3 (1997) 601--642
  \MR{1613528}

\bibitem{SeiD} {\bf P Seidel}, {\it Lectures on four dimensional Dehn twists},
        \arxiv{math.SG/0309012}

\bibitem{TO} {\bf A Tralle}, {\bf J Oprea}, {\it Symplectic manifolds with no 
       K\"ahler structure}, Springer Lecture Notes series 1661 (1997)
  \MR{1465676}
\end{thebibliography}
\end{document}